\documentclass[MSNbibl,citesort,seceqn,number,dvips]{arxbj}

\aid{0}
\volume{19}
\issue{3}
\pubyear{2013}
\firstpage{721}
\lastpage{747}
\doi{10.3150/12-BEJ464} 

\makeatletter
\newcommand{\rrvert}{\vert}
\newcommand{\llvert}{\vert}
\newremark{remark}{Remark}

\newtheorem{theorem}{Theorem}[section]
\newtheorem{lemma}[theorem]{Lemma}
\newtheorem{lem}{Lemma}
\newtheorem{proposition}[theorem]{Proposition}
\newproclaim{assum}{Assumption}
\newproclaim{cond}{Condition}

\newcommand{\trace}{\operatorname{trace}}
\renewcommand{\mid}{\vert}
\makeatother

\begin{document}
\begin{frontmatter}

\title{Single index regression models in the presence of censoring
depending on the~covariates}
\runtitle{Censored single index regression}

\begin{aug}
\author[1]{\fnms{Olivier} \snm{Lopez}\corref{}\thanksref{1}\ead[label=e1]{olivier.lopez0@upmc.fr}},
\author[2]{\fnms{Valentin} \snm{Patilea}\thanksref{2}\ead[label=e2]{patilea@ensai.fr}} \and
\author[3]{\fnms{Ingrid} \snm{Van Keilegom}\thanksref{3}\ead[label=e3]{ingrid.vankeilegom@uclouvain.be}}
\runauthor{O. Lopez, V. Patilea and I. Van Keilegom} 
\address[1]{Universit\'e Pierre et Marie Curie, Laboratoire de
Statistique Th\'eorique et Appliqu\'ee, 4 place Jussieu Barre 15-25
Etage 2,
75005 Paris, France. \printead{e1}}
\address[2]{CREST (Ensai) \& IRMAR, Campus de Ker-Lann, Rue Blaise
Pascal BP 37203, 35172 Bruz cedex, France. \printead{e2}}
\address[3]{Institute of Statistics, Universit\'e catholique de Louvain,
Voie du Roman Pays 20, 1348 Louvain-la-Neuve, Belgium. \printead{e3}}
\end{aug}

\received{\smonth{2} \syear{2010}}
\revised{\smonth{7} \syear{2012}}

%
\begin{abstract}
Consider a random vector $(X',Y)'$, where $X$ is $d$-dimensional and
$Y$ is one-dimensional.
We assume that $Y$ is subject to random right censoring. The aim of
this paper is twofold.
First, we propose a new estimator of the joint distribution of
$(X',Y)'$. This estimator overcomes the common
curse-of-dimensionality problem, by using a new dimension reduction
technique. Second, we assume that the relation
between $X$ and $Y$ is given by a {mean regression} single index
model, and propose a new estimator of the parameters in this model.
The asymptotic properties of all proposed estimators are obtained.
\end{abstract}

%

\end{frontmatter}

\section{Introduction and model}

Consider a random vector $(X',Y)'$, where $X=(X^{(1)},\ldots,X^{(d)})'$
is $d$-dimensional and $Y$ is one-dimensional. We assume that $Y$ is
subject to random right censoring, that is, instead of observing
$(X',Y)'$, we observe the triplet $(X',T,\delta)'$, where $T = Y
\wedge
C$, $\delta=\mathbf{1}_{Y \le C}$, and the random variable $C$ is the
censoring variable. Typically, $Y$ is (a transformation of) the
survival time (whose range can span the whole real line), and $X$ is a
vector of characteristics. The data consist of $n$ i.i.d. replications
$(X_i',T_i,\delta_i)'$ of $(X',T,\delta)'$.

Under this setting, the purpose of this paper is twofold. First, we
propose a new
estimator of the joint distribution $F(x,y) = \mathbb{P}(X \le x, Y
\le
y)$ of $X$ and $Y$ (where $X \le x$ means that { $X^{(j)} \leq
x^{(j)}$} for $j=1,\ldots,d$). Second, we assume that the
relation between $X$ and $Y$ is given by a single index mean regression
model (as in, e.g., H\"ardle and Stoker~\cite{r12}, Powell, Stock and
Stoker~\cite{r25}, Ichimura~\cite{r16}, H\"ardle, Hall and Ichimura~\cite{r11},
Klein and Spady~\cite{r17}, Horowitz and H\"ardle~\cite{r14}, Hristache,
Juditsky and Spokoiny~\cite{r15}), and we propose new estimators of 
the\vadjust{\goodbreak} parameters under
this model. These estimators will be constructed under the following
fundamental model assumption on the relation between $Y$ and $C$, which
we impose throughout this paper:

\begin{enumerate}[(A0)]
\item[(A0)] There exists a function $g \dvtx  \mathbb{R}^d
\rightarrow
\mathbb{R}$, such that:
\begin{enumerate}[(ii)]
\item[(i)] $Y$ and $C$ are independent, conditionally on $g(X)$
\item[(ii)] $\mathbb{P}(Y \le C|X,Y) = \mathbb{P}(Y \le C|g(X),Y)$.
\end{enumerate}
\end{enumerate}

Note that assumption (A0) holds in the particular case
where $\mathcal{L}(C|X,Y)=\mathcal{L}(C|g(X)).$ By assuming that
the censoring variable depends on $X$ only through a
one-dimensional variable $g(X)$, we avoid the curse-of-dimensionality problems
which strike regression approaches where $X$ is multivariate and
$Y$ is independent of $C$ conditionally on $X$, and at the same time
the dependence of $C$ on $X$ is not too restrictive. A related
dimension reduction model assumption for the censoring time has been
considered in Section 4 of Li, Wang and Chen~\cite{r18}.

The function $g$ will be unknown in general. When $g$ is known,
assumption (A0) has been proposed by Lopez~\cite{r19}.
The assumption is needed for identifying the model. In the literature on
nonparametric censored regression, alternatives to assumption (A0)
have been proposed. There are basically two alternatives, which can be
regarded as limiting cases of assumption (A0), and in that sense our
assumption is a trade-off between these two. The first alternative has
been used by, for example, Akritas~\cite{r1} and Van Keilegom and Akritas~\cite{r33}, among many others.
They assume that
$Y$ is independent of~$C$, conditionally on $X$, and propose
kernel type estimators of the distribution $F(x,y)$ under this
assumption. This assumption is a particular case of (A0) by taking
$g(X) \equiv X$.
Their estimators are however restricted to the case where $d=1$.
Although they
could in principle be extended to higher dimensions, this is not
recommended in practice, since they will suffer from the curse-of-dimensionality
and higher order kernels will need to be used. The second alternative to
assumption (A0) has been proposed by Stute~\cite{r28,r29}. He assumes that
$Y$ is independent of $C$, and that $\mathbb{P}(Y\leq
C|X,Y)=\mathbb{P}(Y\leq C|Y)$. This is again a particular case of (A0),
by taking $g(X) \equiv1$.
Although his estimator can be used for any $d \ge1$, it has the
drawback that
it assumes that the censoring variable $C$ depends on $X$ in a very
particular way.
This type of dependence might hold true when the censoring is purely
`administrative'
(censoring at the end of the study), but when the censoring can be
caused by other factors
(like death due to another disease, change of treatment, \ldots), then
less restrictive
assumptions on the censoring mechanism are required.

Our assumption (A0) balances
somewhere in between these two extreme assumptions.
By imposing assumption (A0), we propose a new dimension reduction technique,
which overcomes the drawbacks of these two classical sets of assumptions,
by allowing for $d \ge1$ without assuming the complete independence
between $Y$ and $C$.

In some cases, the function $g$ will be known exactly from some
a priori information. For example, we might know that the censoring
only depends on one component of $X$, for example, $g(X)=X^{(1)}$.
Lopez~\cite{r19} proposed an estimator of the joint distribution $F(x,y)$
when $g$ is supposed to be known.
However, in many other cases, $g$ will be unknown and needs to be estimated.
Throughout this paper, we will assume that
%
\begin{eqnarray}
\label{classG} g \in\mathcal{G}, \qquad \mbox{where } \mathcal{G} = \bigl\{x
\rightarrow\lambda (\theta,x) \dvt \theta\in\Theta \bigr\},
\end{eqnarray}
where $\lambda$ is a known function, and $\Theta$ is a compact
parameter set in
$\mathbb{R}^k$. The true (but unknown) value of $\theta$ will be
denoted by $\theta_0$.
This semiparametric assumption on the conditional distribution of $C$ allows
to avoid the curse of dimensionality that would have stroke our
approach if no restriction on the
censoring time would have been made.

Throughout the paper, we will assume that we know some root-$n$
consistent estimator
$\hat{\theta}$ of $\theta_0,$ that satisfies the following:
\begin{enumerate}[(C0)]
\item[(C0)] The estimator $\hat\theta$ satisfies:
\[
\hat{\theta} -\theta_0 = \frac{1}{n}\sum
_{i=1}^n \mu(T_i,\delta_i,X_i)+
\mathrm{o}_P \bigl(n^{-1/2} \bigr),
\]
with $E[\mu(T,\delta,X)]=0$ and $E[\mu(T,\delta,X)^2]<\infty.$
\end{enumerate}
Hence, the set $\Theta$ {can from now on be an arbitrarily} small
environment of $\theta_0$.

To illustrate the nature of assumptions (A0) and (C0), consider the function
$g(x)=\theta_0'x,$ and the case where $C$ follows a Cox regression
model given $X$, in the sense that the conditional
hazard $h(\cdot|x,y)$ of $C$ given $X=x$ and $Y=y$ satisfies
\[
h(c|x,y) = h_0(c) \exp \bigl(\theta_0'x
\bigr)
\]
for some baseline function $h_0$ only depending on $c$.
{This model assumption on $C$ seems realistic}
since often the censoring variable $C$ represents itself a lifetime,
like the time until a patient dies from a disease other than the
disease under study. Under this model, we
clearly have $\mathcal{L}(C|X,Y)=\mathcal{L}(C|\theta_0'X),$ and
the estimator $\hat{\theta}$ proposed by Andersen and Gill~\cite{r3}
satisfies condition (C0), with
\[
\mu(t,\delta,x) = \Sigma^{-1} \biggl((1-\delta)\phi(x,t)-\int
\phi(x,u)\mathbf {1}_{t>u}{ \bigl[1-G(u-|x) \bigr]^{-1}\,
\mathrm{d}G(u|x)} \biggr),
\]
where the matrix $\Sigma$ is defined by condition D in Andersen and Gill~\cite{r3},
\[
\phi(x,t) = x-{\frac{E[X\mathrm{e}^{\theta_0'X}(1-H(t|X))]}{E[\mathrm{e}^{\theta_0'X}(1-H(t|X))]}},
\]
with $H(t|x) = \mathbb{P}(T \le t|X=x)$ and ${G(c|x) =
\mathbb{P}(C \leq c|X=x)}$. See also Gorgens and Horowitz~\cite{r10} for
regression models more general than Cox in which
$\mathcal{L}(C|X,Y)=\mathcal{L}(C|\theta_0'X)$. Alternatively, one
could also assume that { $C = r(\theta_0'X) + U$, where
$r(\cdot)$ is given,} $E(U)=0$, and $U$ is independent of $X$ and
$Y$. For the estimation of $\theta_0$ and the verification of
condition (C0) under this model, see, for example, Akritas and Van Keilegom~\cite{r2} and Heuchenne and Van Keilegom~\cite{r13}.

The purpose of this paper is twofold. The first contribution of
this paper consists in proposing and studying a new nonparametric
estimator of the joint distribution of $X$ and $Y$ under assumption
(A0). Under different sets of assumptions on the relation between
$X$, $Y$ and $C$, this distribution has been the object of study of
many papers in the past. See, for example, Akritas~\cite{r1}, Stute~\cite{r28,r29}, Van Keilegom and Akritas~\cite{r33}, among others.
As mentioned before, assumption (A0)\vadjust{\goodbreak}
allows to avoid the curse-of-dimensionality problem present in some
of these contributions, and the heavy assumptions on the relation
between $C$ and $X$, which are present in many others.

The second contribution of this paper is the estimation of a
semiparametric single index regression model for the
censored response $Y$ given $X$ under assumption (A0). The proposed
estimator is based on a two-step procedure, in which first a
preliminary (consistent) estimator is obtained, which is then used
to build a least squares criterion that defines our new semiparametric
estimator in order to achieve $n^{1/2}$-consistency. Both steps
heavily rely on the estimator of $F(x,y)$ studied before. Note that
in this second contribution two dimension reduction techniques are
used: the first one comes from assumption (A0), which is concerned
with the relation between $Y$ and $C$, and the second one comes from
the single index model, which is making an hypothesis on the
relation between $Y$ and $X$.

Single index regression models are now a common
semiparametric multivariate explanatory approach, see for instance
Delecroix, Hristache and Patilea~\cite{r5} for a review. However, the
literature on single index models with a censored response variable
is rather poor. To the best of our knowledge, the only contribution
that allows for a general relationship between the censoring
variable and the covariates is Li, Wang and Chen~\cite{r18} and it is
based on sliced inverse regression (SIR). However, it is well known
that the SIR approach requires a linear conditional expectation
condition among the covariates, which {may be restrictive in
applications}, see equation (2.3) in Li, Wang and Chen~\cite{r18}.

Lopez~\cite{r20} proposed a semiparametric least squares estimator for
the single index regression in the particular case where $g(X)\equiv
1$ in assumption (A0). A similar procedure was introduced by Wang \emph{et al.}~\cite{r34} under the stronger assumption that $C$ is
independent of $(X^\prime, Y)^\prime$. See also Lu and Cheng~\cite{r23}.
Lu and Burke~\cite{r22} used the same more restrictive condition to
define an average derivative estimator of the index. It is
worthwhile to notice that these three contributions involve a
Kaplan--Meier estimate of the censoring distribution, while in
general assumption (A0) requires a nonparametric estimate of the
\emph{conditional} distribution of $C$ given $g(X)$.

This paper is organized as follows. In the next section, the
estimators of the joint distribution and of the parameters in the
single index model are explained in detail. Section~\ref{sec3} is devoted to
the presentation of the asymptotic results of the proposed
estimators, while in Section~\ref{sec4} we compare our estimator with an
existing estimator in the literature. Finally, Appendix~\ref{appA} contains the
assumptions under which
the results of Section~\ref{sec3} are valid, while Appendix~\ref{appB} contains some
technical lemmas and the proofs of the main results.\looseness=-1

\section{The estimators}

\subsection{Estimation of the distribution $F(x,y)$}

We first explain how to estimate the joint distribution $F(x,y)$ of $X$
and $Y$.
For an arbitrary value of $\theta$, let
%
\begin{equation}
\label{gtheta} G_{\theta}(t|z)=\mathbb{P} \bigl(C\leq t|\lambda(
\theta,X)=z \bigr),
\end{equation}
and define
%
\begin{equation}
\label{beran1} \hat{G}_{\theta}(t|z)=1-\prod_{T_i\leq t}
\biggl(1-\frac{w^{\theta
}_{in}(z)}{\sum_{j=1}^nw^{\theta}_{jn}(z)\mathbf{1}_{T_j\geq
T_i}} \biggr)^{1-\delta_i},
\end{equation}
where
\[
w^{\theta}_{in}(z)={K \biggl(\frac{\lambda(\theta
,X_i)-z}{a_n} \biggr)}\bigg/{\sum
_{j=1}^nK \biggl(\frac{\lambda(\theta,X_j)-z}{a_n}
\biggr)}.
\]
Here, $a_n$ is a bandwidth sequence converging to zero as $n$ tends to
infinity, and $K$ is a probability density function (kernel). Note that
$\hat G_\theta(t|z)$ reduces to the estimator proposed by Beran~\cite{r4}
when $\lambda(\theta,X)$ is equal to $X$.

With at hand the estimator $\hat{\theta}$ introduced in condition (C0),
and the corresponding estimator $\hat{g}(x)=\lambda(\hat{\theta
},x)$ of $g(x),$
we now define the following estimator of $F(x,y)$:
%
\begin{eqnarray}
\hat{F}_{\hat{g}}(x,y) &=& \frac{1}{n}\sum
_{i=1}^n \frac{\delta_i\mathbf{1}_{T_i\leq y, X_i\leq
x}}{1-\hat{G}_{\hat{\theta}}(T_i-|\hat{g}(X_i))}. \label{estimatedg}
\end{eqnarray}
Note that this estimator is in the same spirit as the estimator proposed
by Stute~\cite{r28,r29}, but the denominators of the two estimators
are different, because of the different sets of underlying assumptions.
See also Fan and Gijbels~\cite{r9} for a similar weighting scheme in a
nonparametric regression framework.
Also note that when $g$ would be known, this estimator equals the estimator
proposed and studied in Lopez~\cite{r21}.

In Section~\ref{sec3.1}, we will study the asymptotic properties of the estimator
$\hat{F}_{\hat{g}}(x,y)$.

\subsection{Estimation of the single index model}

We first need to introduce some notations. For $\theta\in\Theta$, let
$Z_{\theta}=\lambda(\theta,X)$,
and let $\mathcal{Z}_
\theta\subset\mathbb{R}$ be the support of the variable $Z_\theta$.
We assume that $\mathcal{Z}_\theta$ is compact for all $\theta\in
\Theta$.
Also, define $H_{\theta}(t|z) = \Bbb{P}(T\leq t | Z_{\theta}=z)$ and
let $\tau_{H_{\theta},z}=\inf\{t \dvt  H_{\theta}(t|z)=1\}.$

We assume that the following single index mean regression model is valid:
for some $\beta_0\in\mathcal{B} \subset\mathbb{R}^d,$ with, say,
first component $\beta_0^{(1)}=1$,
%
\begin{equation}
\label{sim} E[Y\mid X, Y\leq\tau]=E \bigl[Y\mid\beta_0'X,
Y\leq\tau \bigr]=m \bigl(\beta_0'X \bigr),
\end{equation}
where $m$ is an unknown function, and where $\tau$ is some fixed
truncation point, satisfying
\[
\tau<\inf_{\theta\in\Theta} \inf_{z\in\mathcal{Z}_\theta} \tau_{H_{\theta},z}.
\]
Let $f(t;\beta)=E[Y|\beta'X=t,Y\leq\tau].$ Then,
$f(\cdot;\beta_0)=m(\cdot)$. Also, let $\mathcal{B}=\{1\}\times
\tilde{\mathcal{B}},$ where $\tilde{\mathcal{B}}$ is a compact
subset of $\mathbb{R}^{d-1}$, and denote by $\mathcal{X}$ the
support of the covariate vector $X$, which is a compact subset of
$\mathbb{R}^d$.

The truncation at $\tau$ in model (\ref{sim}) is very common in the
context of regression with right censored observations, and is caused
by the lack of information in the right tail of the conditional
distribution of $Y$ given $X$.
See, for example, Akritas~\cite{r1} and Akritas and Van Keilegom~\cite{r2} for
similar truncation mechanisms. Note that when $\mathcal
{L}(Y|X)=\mathcal
{L}(Y|\beta_0'X)$, that is, when the whole distribution of $Y$ given
$X$ only depends on $X$ via $\beta_0'X$, then model (\ref{sim}) is
satisfied for any value of
$\tau$.

The estimation of $\beta_0$ consists of several steps. We first explain
these steps in an informal, intuitive way to outline the main ideas
behind the proposed method, and we next work out each of these steps in
a rigorous way.
\begin{enumerate}
\item Estimate $f(t;\beta)$ using some
nonparametric estimator $\hat{f}(t;\beta)$.
\item Construct a preliminary consistent estimator $\beta_n$ of $\beta_0$.
\item Use $\beta_n$ to compute a {trimming function that helps to avoid
technical
problems caused by denominators close to zero} in the nonparametric
estimation of $f(t;\beta)$.
\item Construct a second semi-parametric estimator $\hat\beta$ of
$\beta_0$
by using the trimming function of the preceding step.
\end{enumerate}

\subsubsection{\texorpdfstring{Estimation of $f(t;\beta)$}{Estimation of f(t;beta)}}

One possible estimator of $f(t;\beta)$ is
%
\begin{equation}
\label{estim} \hat{f}(t;\beta) = {\int \tilde K\biggl(\frac{\beta'x-t}{h}\biggr)y
\mathbf{1}_{y\leq
\tau}\,\mathrm{d}\hat{F}_{\hat{g}}(x,y)}\bigg/\biggl({\int
\tilde K\biggl(\frac{\beta'x-t}{h}\biggr)\mathbf{1}_{y\leq
\tau}\,\mathrm{d}
\hat{F}_{\hat{g}}(x,y)}\biggr),
\end{equation}
where  $h=h_n$ is a second bandwidth sequence,
possibly different from the bandwidth $a_n$ used to estimate the
joint distribution $F(x,y)$, and where $\tilde K$ is a kernel
function. However, other estimators may be used, for
example, $[\hat{F}_{\beta}(\tau|t)]^{-1} \int y \mathbf{1}_{y\leq
\tau} \,\mathrm{d}\hat{F}_{\beta}(y|t)$, where $\hat{F}_{\beta}(y|t)$ denotes
Beran's~\cite{r4} estimator of $\mathbb{P}(Y\leq y|\beta'X=t)$.

In what follows, we do not specify the choice of estimator of
$f(t;\beta)$.
Instead we will work with a generic estimator $\hat f(t;\beta)$ that
satisfies certain
conditions that need to be fulfilled in order to obtain the asymptotic normality
of $\hat\beta$, and we will prove in Section~\ref{sec3.2} that the estimator in
(\ref{estim})
satisfies these conditions.


\subsubsection{\texorpdfstring{Preliminary estimation of $\beta_0$}{Preliminary estimation of beta 0}}

We assume that we know some set $B$ such that
\[
\inf_{\beta\in\mathcal{B}, x\in B}f_\beta^\tau \bigl(\beta'x
\bigr)=c>0,
\]
where the function $f_\beta^\tau$ denotes the density of $\beta'X$,
conditionally
on $Y \le\tau$. Define the following preliminary trimming function:
%
\begin{equation}
\label{trimmingpreliminaire} \tilde{J}(x)=\mathbf{1}_{x\in B}.
\end{equation}
%
Let $M(\beta,f,\tilde{J})=E[(Y-f(\beta'X;\beta))^2\mathbf
{1}_{Y\leq\tau
}\tilde{J}(X)],$
and note that this is minimized as a function of $\beta$ when $\beta
=\beta_0$. Motivated by this fact, we
define the preliminary estimator $\beta_n$ of $\beta_0$ by replacing
all unknown\vadjust{\goodbreak} quantities in $M(\beta,f,\tilde{J})$ by appropriate
estimators, that is,
%
\begin{eqnarray}
\label{estimpreliminaire} \beta_n &=& \arg\min_{\beta\in\mathcal{B}} \int
\bigl(y-\hat{f} \bigl(\beta'x;\beta \bigr) \bigr)^2
\mathbf{1}_{y\leq\tau} \tilde {J}(x)\,\mathrm{d}\hat {F}_{\hat{g}}(x,y)
\nonumber
\\[-9.5pt]
\\[-9.5pt]
&=& \arg\min_{\beta\in\mathcal{B}} M_n(\beta,\hat{f},\tilde{J}).
\nonumber
\end{eqnarray}
Note that other criterion functions can be used, based on $M$ or
$L$-estimating functions. We do not consider them here, since their
analysis is very similar to the one for the least squares criterion
function.\vspace*{-3pt}

\subsubsection{New trimming function}\vspace*{-2pt}

We will now refine the definition of the trimming function, by using
the preliminary estimator $\beta_n$. Define
%
\begin{equation}
\label{nouveautrimming} J(x)=\mathbf{1}_{f^{\tau}_{\beta_n}(\beta_n'x)>c},
\end{equation}
so instead of requiring that $f_\beta^\tau(\beta'x)>c$ for all
$\beta$,
we now only consider $\beta=\beta_n$,\vspace*{1pt} which will be satisfied for many
more $x$-values, and hence this new function $J(x)$ is trimming much
less than the preliminary naive trimming function $\tilde J(x)$.

To simplify our discussion, we will directly consider that the
true function $f^{\tau}_{\beta_n}$ is used in the definition of
$J$. In practice, the trimming function can be estimated by
$\mathbf{1}_{\hat{f}^{\tau}_{\beta_n}(\beta_n'x)>c},$ where
\[
\hat{f}^{\tau}_{\beta}(t)=\frac{1}{nb_n \mathbb{P}(Y\leq\tau
)}\sum
_{i=1}^n \frac{\delta_i\mathbf{1}_{T_i\leq\tau}}{1-\hat{G}_{\hat
\theta
}(T_i-|\hat{g}(X_i))}K \biggl(
\frac{\beta'X_i-t}{b_n} \biggr),
\]
and where $b_n \rightarrow0$ is a bandwidth parameter. In
applications, $c_1=c \mathbb{P}(Y\leq\tau)$ can be chosen
arbitrarily small by the statistician. Considering
$f^{\tau}_{\beta_n}$ or $\hat{f}^{\tau}_{\beta_n}$ does not change\vspace*{1pt}
anything asymptotically speaking, see { the arguments in
Delecroix, Hristache and Patilea~\cite{r5}, see also Step 0 in the proof
of Theorem~\ref{th4} below. By similar arguments,} the estimator of
$\beta_0$ obtained with
$\mathbf{1}_{\hat{f}^{\tau}_{\beta_n}(\beta_n'x)>c}$ is
asymptotically equivalent to the `ideal' estimator obtained with the
trimming function
%
\begin{equation}
\label{nouveautrimming2} J_0(x)=\mathbf{1}_{f^{\tau}_{\beta_0}(\beta_0'x)>c},
\end{equation}
as long as $\beta_n$ is a consistent estimator of $\beta_0$.
{ Let us point out that $J_0$ only depends on $\beta_0'x$
and, in view of equation (\ref{rouge}) in the proof of Theorem
\ref{th4}, this property will be essential for achieving
$\sqrt{n}$-asymptotic normality of our estimator $\hat\beta$
defined below.}\vspace*{-3pt}

\subsubsection{\texorpdfstring{Estimation of $\beta_0$}{Estimation of beta 0}}\vspace*{-2pt}

With at hand this new trimming function, we can now define a new semi-parametric
least squares estimator of $\beta_0$:
%
\begin{eqnarray}
\label{hatbeta} \hat{\beta} &=& \arg\min_{\beta\in\mathcal{B}_n} \int \bigl(y-\hat{f}
\bigl( \beta'x;\beta \bigr) \bigr)^2
\mathbf{1}_{y\leq
\tau}J(x)\,\mathrm{d} \hat{F}_{\hat{g}}(x,y)
\nonumber
\\[-9.5pt]
\\[-9.5pt]
&=& \arg\min_{\beta\in\mathcal{B}_n} M_n(\beta,\hat{f},J),
\nonumber\vadjust{\goodbreak}
\end{eqnarray}
where $\mathcal{B}_n$ is a set shrinking to $\{\beta_0\},$ which
is computed from the preliminary step. The proof of the asymptotic
normality of $\hat{\beta}$ will be carried out in two steps.
We will first show that minimizing $M_n(\beta,\hat{f},J)$ is
asymptotically equivalent to minimizing $M_n(\beta,f,J_0)$.
This then brings back the minimization problem to a fully
parametric one.

\section{Asymptotic properties}\label{sec3}

\subsection{Estimation of the distribution $F(x,y)$}\label{sec3.1}

Let us first introduce a few notations. Denote $H(t) = \Bbb{P}(T \le
t)$, $H_{\theta}(t|z) = \Bbb{P}(T \le t |\break Z_{\theta}=z)$,
$H_{\theta,0}(t|z) = \Bbb{P}(T \le t, \delta=0 | Z_{\theta}=z)$, and
$H_{\theta,1}(t|z) = \Bbb{P}(T \le t, \delta=1 | Z_{\theta}=z)$.
For any function $L(u)$, let $\nabla_u L(u)$ (resp., $\nabla^2_{u,u} L(u)$) denote the vector (resp., matrix) of partial
derivatives of order 1 (resp., order 2) of $L(u)$ with respect
to $u$. In particular, denote by $\nabla_{\theta}G_{\theta
}(t|\lambda
(\theta,x))$ the vector of partial derivatives of the function
$G_{\theta}(t|\lambda(\theta,x))$ with respect to all occurrences of
$\theta$. Let us point out that, in general, the vector valued function
$\nabla_{\theta}G_{\theta}(t|\lambda(\theta,x))$ depends on $x,$ and
not only on $\lambda(\theta,x).$ Finally, for any matrix $A$ of
dimensions $k \times\ell$ (where $k,\ell\ge1$) we denote $|A| =
[\trace(A'A)]^{1/2}$.

We further need to introduce two (intermediate) estimators of $F(x,y)$:
%
\begin{eqnarray}
\label{tildeG} \tilde{F}_g(x,y) &=& \frac{1}{n}\sum
_{i=1}^n \frac{\delta_i\mathbf{1}_{T_i\leq y, X_i\leq
x}}{1-G_{\theta_0}(T_i-|g(X_i))},
\\
\label{trueg} \hat{F}_g(x,y) &=& \frac{1}{n}\sum
_{i=1}^n \frac{\delta_i\mathbf{1}_{T_i\leq y, X_i\leq
x}}{1-\hat{G}_{\theta_0}(T_i-|g(X_i))}.
\end{eqnarray}

In the following result, we consider integrals of the form $\int\phi
( x,y )
\,\mathrm{d}\hat{F}_g ( x,y )$ with $\phi$ belonging to some class of
functions $\mathcal{F}$, and we state that this class of integrals is
Glivenko--Cantelli and admits an i.i.d. representation uniformly over
all $\phi\in\mathcal{F}$. The proof can be found in Lopez~\cite{r21}. For
a completely nonparametric estimator of $F(x,y)$ that is not based on
model assumption (A0), S\'anchez-Sellero, Gonz\'alez-Manteiga and Van Keilegom~\cite{r26} obtained a similar uniform consistency and convergence
result. The assumptions mentioned below can be found in Appendix~\ref{appA}.

\begin{theorem}
\label{t1} \emph{(i)} Under Assumptions~\ref{a1} and~\ref{a8},
for $a_n\rightarrow0$ and $na_n\rightarrow\infty$, and for a class
$\mathcal{F}$ satisfying Condition~\ref{a10}, we have
\[
\sup_{\phi\in\mathcal{F}} \biggl\llvert \int\phi ( x,y ) \,\mathrm{d}[
\hat{F}_g-F] ( x,y ) \biggr\rrvert \rightarrow_{a.s.}0.
\]

\emph{(ii)} For $Z_i=\lambda(\theta_0,X_i),$ define
\[
M_i(t)=(1-\delta_i) \mathbf{1}_{T_i\leq
t}-\int
_{-\infty}^{t}\frac{\mathbf{1}_{T_i\geq
y}\,\mathrm{d}G_{\theta_0}(y|Z_i)}{1-G_{\theta_0}(y-|Z_i)},
\]
which is a continuous time martingale with respect to the natural
filtration $\sigma(\{Z_i\mathbf{1}_{T_i\leq
t},\break T_i\mathbf{1}_{T_i\leq t}, \delta_i\mathbf{1}_{T_i\leq t},
i=1,\ldots,n\})$. Under Assumptions~\ref{a1}--\ref{a9} and for a class
$\mathcal{F}$
satisfying Conditions~\ref{a11} and~\ref{a12},
\[
\int\phi ( x,y ) \,\mathrm{d}[\hat{F}_g-\tilde{F}_g] (
x,y ) = \frac{1}{n}\sum_{i=1}^n \int
\frac{\bar{\phi}(Z_i,s)\,\mathrm{d}M_i(s)}{[1-F(s-\mid Z_i)][1-G_{\theta
_0}(s\mid
Z_i)]} +R_{n} ( \phi ) ,
\]
where $\sup_{\phi\in\mathcal{F}}| R_{n}(\phi)| =\mathrm{o}_P(n^{-1/2})$,
$\bar
{\phi}$ is defined above Condition~\ref{a12}, and $F(s|z)=\mathbb
{P}(Y\leq s|Z_{\theta_0}=z).$
\end{theorem}

The following theorem {provides} the behavior of the difference
between integrals with respect to $\hat{F}_{\hat{g}}$ and
integrals with respect to $\hat{F}_g$.

\begin{theorem}
\label{curse3} \emph{(i)} Under Assumptions~\ref{a1},~\ref{a8} and~\ref{a13},
for $a_n\rightarrow0$ and $na_n\rightarrow\infty$, and for a class
$\mathcal{F}$
satisfying Condition~\ref{a10}, we have
\[
\sup_{\phi\in\mathcal{F}} \biggl\llvert \int\phi(x,y)\,\mathrm{d}[
\hat{F}_{\hat
{g}}-\hat {F}_{g}](x,y) \biggr\rrvert =
\mathrm{o}_P(1).
\]

\emph{(ii)} Under Assumptions~\ref{a1}--\ref{a8} and~\ref{a13}, for
$a_n\rightarrow0$ and $na_n^3(\log n)^{-1}\rightarrow\infty$, and for
a class $\mathcal{F}$ whose envelope
is as in Condition~\ref{a10},
\begin{eqnarray*}
&& \int\phi(x,y)\,\mathrm{d}[\hat{F}_{\hat{g}}-\hat{F}_g](x,y)
\\
&&\quad =-E \biggl(\frac{\phi(X,Y)\{\nabla_{\theta}G_{\theta_0}(Y-|\lambda
(\theta_0,X))\}'}{1-G_{\theta_0}(Y-|g(X))} \biggr) \frac{1}{n}\sum
_{i=1}^n \mu(T_i,\delta_i,X_i)
+ \tilde{R}_n(\phi),
\end{eqnarray*}
where the function $\mu$ is defined in \emph{(C0)}, and where $\sup_{\phi
\in
\mathcal{F}}|\tilde{R}_n(\phi)|=\mathrm{o}_P(n^{-1/2}).$
\end{theorem}

\subsection{Estimation of the single index model}\label{sec3.2}

We now return to the single index model (\ref{sim}) and to the
estimators $\beta_n$ and $\hat\beta$ defined in (\ref
{estimpreliminaire}) and (\ref{hatbeta}). We start with stating the
asymptotic consistency of the estimator $\beta_n$. Note that the
estimator $\hat\beta$ is by construction consistent, since it is
defined on a shrinking neighborhood of $\beta_0$.

\begin{theorem} \label{th3}
Let $\tilde{J}$ be defined as in (\ref{trimmingpreliminaire}). Under
Assumptions~\ref{a1},~\ref{a8},~\ref{a13},~\ref{etlemodele}, and
\ref{premierevitesse} -- \emph{(\ref{reg1})}, and for $a_n\rightarrow0$ and
$na_n\rightarrow\infty$, we have
\[
\sup_{\beta\in\mathcal{B}}\bigl|M_n(\beta,\hat{f},\tilde{J})-M(\beta ,f,
\tilde{J})\bigr|\rightarrow0,
\]
in probability. Consequently,
$\beta_n \rightarrow\beta_0$ in probability.\vadjust{\goodbreak}
\end{theorem}

The next lemma is an important property in the literature on single
index models. In the classical uncensored single index regression
model, the property
$E[\nabla_{\beta}f(\beta_0'X;\beta_0)|\break\beta_0'X]=0$ plays
a major role in proving the asymptotic normality of $M$-estimators.
See Delecroix, Hristache and Patilea~\cite{r5}. The next lemma shows
that in our context, where we have to truncate at $\tau$ because of
censoring in the data, the analogous truncated version of this
property holds true { without any further model
conditions.}

\begin{lemma}
\label{esperanceestnulle} Assume that the derivative $\nabla_{\beta
}f(\beta_0'\cdot;\beta_0)$ exists and is bounded.
Then, for any $\beta_0$ satisfying condition (\ref{sim}),
\[
E \bigl[\nabla_{\beta}f \bigl(\beta_0'X;
\beta_0 \bigr)\mathbf {1}_{Y\leq
\tau}\mid\beta_0'X
\bigr]=0.
\]
%
\end{lemma}

This lemma is crucial for obtaining our i.i.d. representation and the
asymptotic normality of $\hat\beta$, which we state in the next theorem.
We denote by $\nabla_{\tilde\beta} f(\beta_0'\cdot;\beta_0)$ the vector
of partial derivatives with respect to the last $d-1$ components of
$\beta$.

\begin{theorem} \label{th4}
\label{normasymlastone} Let
$\phi(x,y)=(y-f(\beta_0'x;\beta_0))\nabla_{\tilde\beta}f(\beta_0'x;\beta_0) \mathbf{1}_{y\leq\tau} J_0(x)$. Under Assumptions
\ref{a1}--\ref{unpeudevitesse?}, we have
%
\begin{eqnarray}
\label{repiid} \hat{\tilde{\beta}}-\tilde\beta_0 &=&
\Omega^{-1} \Biggl[ \int\phi(x,y) \,\mathrm{d} \bigl(
\tilde{F}_g ( x,y )-F(x,y) \bigr)
\nonumber
\\
&& \hspace*{23pt} {} +\frac{1}{n}\sum_{i=1}^n
\int \frac{\bar{\phi}(g(X_i),s)\,\mathrm{d}M_i(s)}{[1-F(s-| g(X_i))][1-G_{\theta
_0}(s|g(X_i))]}
\nonumber
\\
&&\hspace*{23pt} {} - E \biggl(\frac{\phi(X,Y)\{\nabla_{\theta
}G_{\theta_0}(Y-|\lambda(\theta_0,X))\}'}{1-G_{\theta
_0}(Y-|g(X))} \biggr) \frac{1}{n}\sum
_{i=1}^n \mu(T_i,
\delta_i,X_i) \Biggr]
\\
&&{} + \mathrm{o}_P \bigl(n^{-1/2} \bigr)
\nonumber
\\
&=&
\nonumber
\Omega^{-1} \Biggl[\frac{1}{n}\sum
_{i=1}^n \eta(T_i,\delta_i,X_i)
\Biggr]+\mathrm{o}_P \bigl(n^{-1/2} \bigr),
\end{eqnarray}
where the function $\mu$ is defined in \emph{(C0)}, and where
\begin{eqnarray*}
\Omega&=& E \bigl[\mathbf{1}_{Y\leq\tau} J_0(X)
\nabla_{\tilde\beta}f \bigl(\beta_0'X;
\beta_0 \bigr)\nabla_{\tilde\beta
}f \bigl(\beta_0'X;
\beta_0 \bigr)' \bigr].
\end{eqnarray*}
Hence,
\[
n^{1/2} (\hat{\tilde\beta}-\tilde\beta_0) \stackrel {d} {
\rightarrow} N \bigl(0, \Omega^{-1} E \bigl[\eta(T,\delta,X) \eta(T,
\delta,X)' \bigr] \Omega^{-1} \bigr).
\]
\end{theorem}

If we wish to estimate the asymptotic variance in Theorem~\ref{normasymlastone},
we see that we need to estimate the variance of $\Omega^{-1}\eta.$
However, one can
consistently estimate $\Omega$ by
\[
\hat{\Omega}=\frac{1}{n}\sum_{i=1}^n
\mathbf{1}_{Y_i\leq\tau} J(X_i) \nabla_{\tilde\beta}\hat{f}
\bigl( \hat{\beta}'X_i;\hat{\beta} \bigr)
\nabla_{\tilde
\beta} \hat{f} \bigl(\hat{\beta}'X_i;
\hat{\beta} \bigr)'.
\]
Similarly, when it comes to estimate the covariance matrix of $\eta,$
one can proceed by taking the empirical variance of
a random vector $(\hat{\eta}(T_i,\delta_i,X_i))_{1\leq i \leq n},$
where $\hat{\eta}$ denotes an estimated version of
$\eta$ in which we replaced each unknown quantity by its empirical
counterpart ($f$ replaced by $\hat{f},$ $\beta_0$ by $\hat{\beta},$ $F$
by $\hat{F}, \ldots$).

We end this section with the verification of Assumptions~\ref{premierevitesse}--\ref{unpeudevitesse?} for the estimator $\hat
f(t;\beta)$
defined in (\ref{estim}).
Define the (uncomputable) kernel estimator based on $\tilde{F}_g,$
%
\begin{eqnarray}
\label{f*} f^*(t;\beta) = {\int\tilde K\biggl(\frac{\beta'x-t}{h}\biggr)y\mathbf
{1}_{y\leq\tau}\,\mathrm{d}\tilde{F}_g(x,y)}\bigg/\biggl({\int\tilde K
\biggl(\frac{\beta'x-t}{h}\biggr)\mathbf{1}_{y\leq\tau}\, \mathrm{d}
\tilde{F}_g(x,y)}\biggr).\quad
\end{eqnarray}
The advantage of $\tilde{F}_g$, and hence of $f^*$, is that it is
composed of sums of i.i.d. terms. Classical arguments show that $f^*$
satisfies Assumptions~\ref{premierevitesse} to~\ref{unpeudevitesse?}.
This is shown in Proposition~\ref{tt} below. On the other hand,
Proposition~\ref{diff} shows that the difference between $\hat{f}$ and
$f^*$ is sufficiently small so that $\hat{f}$ also satisfies these assumptions.


\begin{proposition}
\label{tt}
Assume that
\begin{longlist}[(iii)]
\item[(i)] $\tilde{K}$ is a symmetric density function with compact support,
and with two continuous derivatives of bounded variation;
\item[(ii)] $f(\cdot;\beta_0) \in\mathcal{H}_1^0$ and $\nabla_{\beta
}f(\beta_0'\cdot;\beta_0)\in\mathcal{H}_2^0$,
with $\mathcal{H}_1^0$ and $\mathcal{H}_2^0$ defined in \emph{(\ref{H10})} and
\emph{(\ref{H20})};
\item[(iii)] $nh^5(\log n)^{-1/2}\rightarrow\infty,$ and $nh^8\rightarrow
0$.
\end{longlist}
Then, $f^*$ satisfies Assumptions~\ref{premierevitesse}--\ref{unpeudevitesse?}.
\end{proposition}

\begin{proposition}
\label{diff}
Under the assumptions of Theorem~\ref{curse3},
we have
\begin{eqnarray*}
\sup_{\beta\in\mathcal{B},x \in\mathcal{X}}\bigl|f^* \bigl(\beta'x;\beta \bigr)-\hat{f} \bigl(
\beta'x;\beta \bigr)\bigr| &=& \mathrm{O}_P \bigl((\log
n)^{1/2}n^{-1/2}a_n^{-1/2} \bigr),
\\
\sup_{\beta\in\mathcal{B},x \in\mathcal{X}}\bigl|\nabla_{\beta
}f^* \bigl(\beta'x; \beta
\bigr)-\nabla_{\beta}\hat{f} \bigl(\beta'x;\beta \bigr)\bigr| &=&
\mathrm{O}_P \bigl((\log n)^{1/2}h^{-1}n^{-1/2}a_n^{-1/2}
\bigr),
\\
\sup_{\beta\in\mathcal{B},x \in\mathcal{X}}\bigl|\nabla^2_{\beta
,\beta
}f^* \bigl(
\beta'x;\beta \bigr)-\nabla^2_{\beta,\beta}\hat{f}
\bigl(\beta'x;\beta \bigr)\bigr| &=& \mathrm{O}_P \bigl((
\log n)^{1/2}h^{-2}n^{-1/2}a_n^{-1/2}
\bigr),
\end{eqnarray*}
where $\hat f$ is the estimator defined in (\ref{estim}). Moreover, {
$\nabla_{\beta}\hat{f}(\beta_0'x;\beta_0)=x\hat{m}_1(\beta_0'x)+\hat
{m}_2(\beta_0'x),$ with, for $j=1,2,$
\begin{eqnarray*}
\sup_{x \in\mathcal{X}}\bigl|\hat{m}_j \bigl(\beta_0'x
\bigr)-m_j^* \bigl(\beta_0'x \bigr)\bigr|&=&
\mathrm{O}_P \bigl((\log n)^{1/2}h^{-1}n^{-1/2}a_n^{-1/2}
\bigr),
\\
\sup_{u \in\beta_0'\mathcal{X}}\bigl|\hat {m}'_j(u)-m_j^{*\prime}(u)\bigr|&=&
\mathrm{O}_P \bigl((\log n)^{1/2}h^{-2}n^{-1/2}a_n^{-1/2}
\bigr),
\end{eqnarray*}
where the functions $m^*_j$ are defined in \emph{(\ref{decompof*})}, and where
$m'$ denotes the derivative
of the univariate function $\beta_0'\mathcal{X}\ni u\rightarrow m(u).$}
\end{proposition}

Note that $\hat{f}'(u;\beta_0)=\hat{m}_1(u)$ (resp., $f^{*\prime}(u;\beta_0)=m_1^*(u)$). {Combine Propositions~\ref{tt} and~\ref{diff} and
deduce that $\hat{f}$ satisfies Assumptions~\ref{premierevitesse}--\ref{unpeudevitesse?} if
$nh^8 \rightarrow0, {na_nh^4(\log n)^{-1} \rightarrow\infty}$ and
$ha_n^{-1/2}(\log n)^{1/2} \rightarrow0$.}
In the case where $a_n=n^{-1/[4-\delta]}$ for some {$\delta\in(0,1),$}
these conditions are satisfied if $nh^{{4({4-\delta})/({3-\delta
})}}(\log n)^{{-({4-\delta})/({3-\delta})}} \rightarrow\infty$ and
$nh^{8-2\delta}(\log n)^{4-\delta} \rightarrow0$.

\section{Simulation study}\label{sec4}

To investigate the small sample behavior of our procedure, we carry out
a small simulation study in which we consider two models. In the first
model, the regression function is given by
\[
m_1 \bigl(\beta_0'x \bigr)=
\beta_0'x-0.5 \bigl(\beta_0'x
\bigr)^2,
\]
and in the second
\[
m_2 \bigl(\beta_0'x \bigr)=\log
\bigl(1+0.5 \beta_0'x \bigr),
\]
with $\beta_0=(1,0.75,0.25,-0.5).$
We consider residuals $\varepsilon=Y-m_{j}(\beta_0'X)$ (for $j=1,2$)
that are Gaussian variables $\mathcal{N}(0,1)$ independent of $X.$
The covariates are composed of $4$ independent components, following an
uniform distribution on $[0,1].$

The censoring variable $C$ follows an exponential distribution with
mean {$\gamma\exp(\theta_0'X)$ conditional on the covariate $X_i$,
where $\theta_0=(-0.1,-0.2,0.1,-0.3),$ and $\gamma$ }is a parameter
that allows us to modify the
average proportion of censored responses. The parameter $\theta_0$ is
estimated by maximizing the Cox pseudo-likelihood, since the regression
model on $C$ is a proportional hazards model.

\begin{table}[b]
\tabcolsep=0pt
\caption{Comparison of the MSE of the proposed estimator $\hat{\beta}$
(columns CKM) with the MSE of the estimator based on Kaplan--Meier
weights (columns KM) for different proportions of censoring}\label{table1}
\begin{tabular*}{\textwidth}{@{\extracolsep{\fill}}lllllll@{}}
\hline
&\multicolumn{6}{l}{Proportion of censoring}\\ [-5pt]
&\multicolumn{6}{l}{\hrulefill}\\
&\multicolumn{2}{l}{15\%}&\multicolumn{2}{l}{30\%}&\multicolumn{2}{l}{50\%}\\ [-5pt]
&\multicolumn{2}{l}{\hrulefill}&\multicolumn{2}{l}{\hrulefill}&\multicolumn{2}{l}{\hrulefill}\\
Regression model & CKM & KM & CKM & KM & CKM & KM \\
\hline
$m_1$ & 1.022 & 1.463 & 1.147 & 1.279 & 1.619 & 1.728 \\
$m_2$ & 0.580 & 1.480 & 1.290 & 1.613 & 1.407 & 1.633 \\
\hline
\end{tabular*}
\end{table}

We consider $10\mbox{,}000$ replications of this simulation scheme for $n=200.$
For each simulated sample $j$,
we compute the resulting estimator $\hat{\beta}^{(j)}$ of $\beta_0$ and
compute $\|\hat{\beta}^{(j)}-\beta_0\|_2^2.$
We then deduce an estimator of the mean squared error (MSE) $E[\|\hat
{\beta}-\beta_0\|_2^2].$
We take $a_n=2$ for the bandwidth involved in Beran's estimator. Since
the procedure is more sensitive to the choice of the second
bandwidth $h,$ we consider a set of bandwidths $h_j=0.5+j0.1,$ for
$j=1,\ldots,10,$ and for each sample, we take the bandwidth
that gives the lowest value of $M_n(\beta,\hat f, J)$ defined in (\ref
{hatbeta}). In Table~\ref{table1}, we compare the MSE of the estimator
that we propose to the MSE of an estimator based on
Kaplan--Meier weights, that is if we replace Beran's estimator in our
approach by a standard Kaplan--Meier estimator. This alternative
estimator is the one defined in Lopez~\cite{r20}. As for our approach, this
estimator puts more weights to the largest uncensored observations
caused by censoring. Nevertheless this alternative procedure is not
adapted to Assumption (A0) that we use herein. {Hence, the estimator of
Lopez~\cite{r20}
is expected to fail in our simulation setting.}

As expected, our estimator based on the conditional Kaplan--Meier
weighting outperforms the estimator of Lopez~\cite{r20} in the different
situations we consider. {It is also natural to observe that the MSE of
our $\hat{\beta}$ increases with the proportion of censoring.}

\begin{appendix}
\section{Assumptions and conditions}\label{appA}

We split the assumptions in three parts, namely those required for the
estimation of $F(x,y)$, the estimation of $\beta_0$, and the estimation
of $f(\cdot;\beta)$.

\subsection*{Assumptions needed for the estimation of $F(x,y)$}

The asymptotic results related to the estimator $\hat F_{\hat g}(x,y)$
will be valid under the following assumptions and conditions.

\begin{assum} \label{a1}
The distribution $\mathbb{P}(Z_\theta\le z)$ has three uniformly
bounded derivatives for $z\in\mathcal{Z}_\theta$ and $\theta\in
\Theta
$, and the
densities $f_{Z_{\theta}} ( z )$ satisfy $\inf_{\theta\in
\Theta}\inf_{z\in\mathcal{Z}_\theta}f_{Z_{\theta}} ( z
) >0$.
\end{assum}

For any function $J(t\mid z)$ we will denote
by $J_{c}(t\mid z)$ the continuous part, and $J_{d}(t\mid
z)=J(t\mid z)-J_{c}(t\mid z).$ Assumption~\ref{a3} below has been introduced
by Du and Akritas~\cite{r7} to obtain their asymptotic i.i.d.
representation of
the conditional Kaplan--Meier estimator.

\begin{assum}\label{a3}
\emph{(i)} Let $L(y|z)$ denote $H_{\theta_0}(y|z)$ or $H_{\theta
_0,0}(y|z).$ Then, $\nabla_z L(y|z)$ and $\nabla^2_{z,z}L(y|z)$ exist, are
continuous with respect to $z,$ and are uniformly bounded as functions
of $(z,y)$.\vspace*{-6pt}
\begin{longlist}[(iii)]
\item[(ii)]
For some positive nondecreasing bounded (on $ [ -\infty
;\tau ] $) functions $L_{1}$, $L_{2}$, $L_{3}$, we have, for
all $z\in\mathcal{Z}_{\theta_0}$,
\begin{eqnarray*}
\bigl\llvert H_{\theta_0c}(t_{1}\mid z)-H_{\theta_0c}(t_{2}
\mid z) \bigr\rrvert &\leq& \bigl\llvert L_{1} ( t_{1} )
-L_{1} ( t_{2} ) \bigr\rrvert ,
\\
\bigl\llvert \nabla_zH_{\theta_0c}(t_{1}\mid z)-
\nabla_{z}H_{\theta_0c} (t_{2}\mid z) \bigr\rrvert &
\leq& \bigl\llvert L_{2} ( t_{1} ) -L_{2} (
t_{2} ) \bigr\rrvert ,
\\
\bigl\llvert \nabla_z H_{\theta_0,0c}(t_{1}\mid z)-
\nabla_zH_{\theta
_0,0c}(t_{2}\mid z) \bigr\rrvert &\leq&
\bigl\llvert L_{3} ( t_{1} ) -L_{3} (
t_{2} ) \bigr\rrvert ,
\end{eqnarray*}
the last two assumptions implying the same kind for $\nabla_z
H_{1c}$.

\item[(iii)] The jumps of $F_{g}(\cdot\mid z)$ and $G_{\theta
_0}(\cdot
\mid z)$ are the same for all $z\in\mathcal{Z}_{\theta_0}$. Let
$(d_{1},d_{2},\ldots)$ be the atoms of $G$.

\item[(iv)] $F_{g}(\cdot\mid z)$ and $G_{\theta_0}(\cdot\mid z)$ have
two derivatives with respect to
$z$, with the first derivatives uniformly bounded (on $ [
-\infty;\tau ] $). The variation of the functions
$\nabla_zF_{g}(\cdot\mid z)$ and $\nabla^2_{z,z}F_{g}(\cdot\mid z)$
on $ [
-\infty;\tau ] $ is bounded by a constant not depending on
$z$.
\item[(v)] For all $d_{i}$, define
\begin{eqnarray*}
s_{i} &=&\sup_{z\in\mathcal{Z}_{\theta_0}} \bigl\llvert F_g(d_{i}-
\mid z)-F_g(d_{i}\mid z) \bigr\rrvert ,
\\
s_{i}^{\prime} &=&\sup_{z\in\mathcal{Z}_{\theta_0}} \bigl\llvert
\nabla_z F_g(d_{i}-\mid z)-
\nabla_z F_g(d_{i}\mid z) \bigr\rrvert ,
\\
r_{i} &=&\sup_{z\in\mathcal{Z}_{\theta_0}} \bigl\llvert G_{\theta
_0}(d_{i}-
\mid z)-G_{\theta_0}(d_{i}\mid z) \bigr\rrvert ,
\\
r_{i}^{\prime} &=&\sup_{z\in\mathcal{Z}_{\theta_0}} \bigl\llvert
\nabla_z G_{\theta_0}(d_{i}-\mid z)-
\nabla_z G_{\theta_0}(d_{i}\mid z) \bigr\rrvert .
\end{eqnarray*}
Then, $\sum_{d_{i}\leq\tau} (s_{i}+s_{i}^{\prime
}+r_{i}+r_{i}^{\prime
}) <\infty$.
\end{longlist}
\end{assum}
%

\begin{assum} \label{a8} The kernel $K$ is a symmetric probability
density function
with compact support, and $K$ has bounded second derivative.
\end{assum}

\begin{assum} \label{a9} The bandwidth $a_n$ satisfies { $(\log
n) n^{-1}a_n^{-3} \rightarrow0$ } and $na_n^4\rightarrow0.$
\end{assum}

\begin{assum} \label{a13} The function $(x,t,\theta)\mapsto G_{\theta
}(t|\lambda(\theta,x))$ is differentiable with respect to $\theta$, and
the vector $\nabla_{\theta}G_{\theta}(t|\lambda(\theta,x))$ is
uniformly bounded in $(x,t,\theta).$
\end{assum}

The class of functions $\mathcal{F}$ considered in Section~\ref{sec3.1} should
satisfy the following conditions, which are taken over from Lopez~\cite{r21}. The conditions make use of concepts from the context of
empirical processes, which can be found, for example, in Van der Vaart and Wellner~\cite{r32}.

\begin{cond}
\label{a10} Let $p_0(x,y,c)=\mathbf{1}_{y\leq c}[1-G_{\theta
_0}(y-|g(x))]^{-1}.$ The class $p_0\mathcal{F}$ is $\mathbb
{P}_{(X,Y,C)}$-Glivenko--Cantelli, and has an integrable envelope $\Phi_0$ satisfying $\Phi_0(x,y,c)=0$ for $y > \tau$.
\end{cond}

\begin{cond}
\label{a11} The covering number $N(\varepsilon,\mathcal
{F},L^2(\mathbb
{P}_{(X,Y)}))$ is
bounded by $A\varepsilon^{-V}$ for $\varepsilon>0$ and for some $A,
V>0$, and $\mathcal{F}$
has a square integrable envelope $\Phi$ satisfying $\Phi(x,y)=0$ for $y
> \tau$.
\end{cond}

Let $Z=Z_{\theta_0}=g(X)$, let $F_{z}(x,y)=\mathbb{P}(X\leq x,
Y\leq y | Z=z),$ and for any function $\phi(x,y)$, define $\bar{\phi
}(z,s)=\int\mathbf{1}_{s\leq y}\phi(x,y)\,\mathrm{d}F_z(x,y)$.
Let $\mathcal{Z}_{\theta_0,\eta}$ be the set of all points at a distance
at least $\eta>0$ from the complementary of $\mathcal{Z}_{\theta_0}.$

\begin{cond}
\label{a12} 
For
all $\phi\in\mathcal{F},$
$\bar{\phi}$ is
twice differentiable with respect to $z,$ and
\[
\sup_{s\leq
\tau,z\in\mathcal{Z}_{\theta_0,\eta}} \bigl\{\bigl|\nabla_z \bar{\phi }(z,s)\bigr|+\bigl|
\nabla_{z,z}^2 \bar{\phi}(z,s)\bigr| \bigr\}\leq M<\infty
\]
for some constant $M$ not
depending on $\phi$.
Moreover, $\bar{\Phi}$ is bounded on
$\mathcal{Z}_{\theta_0,\eta}\,\times\,]$$-$$\infty;\tau],$ and has bounded
partial derivatives with respect to $z,$ where $\Phi$ is the envelope
function of Condition~\ref{a11}.
\end{cond}

The reason for introducing the set $\mathcal{Z}_{\theta_0,\eta}$ is to
prevent us from boundary effects coming from kernel estimators. See
Lopez~\cite{r21} for a detailed discussion on this issue.

\subsection*{Assumptions needed for the estimation of $\beta_0$}
We next state the additional assumptions needed for the asymptotic
results concerning the estimation of the parameters in the single index model.

\begin{assum}
\label{trimberan}
There exist $0<c_0<c_1<\infty$ and $\eta>0$ such that, for each $c\in
[c_0,c_1]$ and $x \in\mathcal{X}$,
\[
\mathbf{1}_{f^{\tau}_{\beta_0}(\beta_0'x)>c}=1 \quad\Longrightarrow \quad g(x)\in
\mathcal{Z}_{\theta_0,\eta}.
\]
Moreover, assume that
\[
\bigl|f^{\tau}_{\beta_1} \bigl(\beta_1'x
\bigr)-f^{\tau}_{\beta_2} \bigl(\beta_2'x
\bigr)\bigr|\leq C\| \beta_1-\beta_2\|^{\alpha}
\]
for some positive constant $C$ and some $\alpha>0$.
\end{assum}


\begin{assum}
\label{etlemodele}
\emph{(i)}  $E(|Y|^3)<\infty$;\vspace*{-6pt}
\begin{longlist}[(iii)]
\item[(ii)] $E [\{f(\beta'X;\beta)-f(\beta_0'X;\beta_0)\}^2\mathbf
{1}_{Y\leq
\tau} ]=0 \Longrightarrow\beta=\beta_0$;
\item[(iii)] $\beta_0=(1,\tilde{\beta}_0')'$ with $\tilde{\beta}_0$ an
interior point of $\tilde{\mathcal{B}}$;
\item[(iv)] The class $\{(x,y)\rightarrow f(\beta'x;\beta)\mathbf
{1}_{y\leq
\tau} \dvt \beta
\in\mathcal{B}\}$ satisfies Condition~\ref{a10} for a continuous
integrable envelope $\Psi$.
\end{longlist}
\end{assum}


\begin{assum}
\label{encoremodel}
The classes $\{x\rightarrow\nabla_{\beta}f(\beta'x;\beta) \dvt  \beta
\in
\mathcal{B} \}$ and
$\{x\rightarrow\nabla^2_{\beta,\beta}f(\beta'x;\beta) \dvt  \beta\in
\mathcal{B} \}$ are VC-classes of continuous functions for a uniformly
bounded envelope.
\end{assum}

\subsection*{Assumptions needed for the estimation of $f(\cdot;\beta)$}
The last group of assumptions is required for the
generic estimator $\hat f(\cdot;\beta)$. They are verified in Section
\ref{sec3.2} for the estimator defined in (\ref{estim}).

\begin{assum}
\label{premierevitesse} For all $c>0$,
%
\begin{eqnarray}
\label{reg1} \sup_{\beta\in\mathcal{B}, x\in
\mathcal{X}}\bigl|\hat{f} \bigl(\beta'x;\beta
\bigr)-f \bigl( \beta'x;\beta \bigr)\bigr|\mathbf {1}_{f^{\tau
}_{\beta}(\beta'x)>c}&=&
\mathrm{o}_P(1),
\\
\label{derreg1} \sup_{\beta\in\mathcal{B}, x\in\mathcal{X}}\bigl|\nabla_{\beta
}\hat {f} \bigl(
\beta'x; \beta \bigr)-\nabla_{\beta}f \bigl(
\beta'x;\beta \bigr)\bigr|\mathbf {1}_{f^{\tau
}_{\beta}(\beta'x)>c} &=&
\mathrm{o}_P(1),
\\
\label{der2reg1} \sup_{\beta\in\mathcal{B}, x\in\mathcal{X}}\bigl|\nabla^2_{\beta
,\beta
}
\hat{f} \bigl( \beta'x;\beta \bigr)-\nabla^2_{\beta,\beta}f
\bigl( \beta'x;\beta \bigr)\bigr|\mathbf {1}_{f^{\tau}_{\beta}(\beta'x)>c} &=&
\mathrm{o}_P(1).
\end{eqnarray}
\end{assum}




\begin{assum}
\label{classededonsker} There exist Donsker classes $\mathcal{H}_1$ and
$\mathcal{H}_2$ such that $f(\cdot;\beta_0) \in\mathcal{H}_1$ and
$\nabla_{\beta}f(\beta_0'\cdot;\beta_0)\in\mathcal{H}_2$, and such
that with
probability tending to one, $\hat{f}(\cdot;\beta_0) \in\mathcal
{H}_1$ and
$\nabla_{\beta}\hat{f}(\beta_0'\cdot; \beta_0)\in\mathcal{H}_2$.
\end{assum}

Typical examples of such kind of Donsker classes are classes of regular
functions. Let $\mathcal{T} = \{\beta_0'x \dvt  x \in\mathcal{X}\}
\subset
\mathbb{R}$ and let
$\mathcal{C}^1_\ell(\mathcal{T},M) = \{h \dvtx  \mathcal{T} \mapsto
\mathbb
{R}^\ell\dvt  \sup_{t \in \mathcal{T}}\{|h(t)|+|h'(t)|\} \le M\}$ for $\ell
\ge1$ and for some $M < \infty$.
Define
%
\begin{eqnarray}
\label {H10} \mathcal{H}_1^0 &=&\mathcal{C}^1_1(
\mathcal{T},M),
\\
\label{H20} \mathcal{H}_2^0 &=& \bigl\{h \dvtx
\mathcal{X} \mapsto \mathbb{R}^d \dvt x \mapsto x h_1
\bigl(\beta_0'x \bigr) + h_2 \bigl(
\beta_0'x \bigr) \dvt h_1 \in\mathcal
{C}^1_1(\mathcal{T},M), h_2 \in
\mathcal{C}^1_d(\mathcal{T},M) \bigr\}.\quad
\end{eqnarray}
%
The class $\mathcal{H}_2^0$ is a Donsker class, which follows from
stability properties of Donsker classes (see, e.g., Examples 2.10.7
and 2.10.10 in Van der Vaart and Wellner~\cite{r32}).

\begin{assum}
\label{unpeudevitesse?} For all $c>0$,
\begin{eqnarray*}
\sup_{x\in\mathcal{X}} \bigl|\hat{f} \bigl(\beta_0'x;
\beta_0 \bigr)-f \bigl(\beta_0'x;
\beta_0 \bigr)\bigr|\mathbf{1}_{f^{\tau}_{\beta_0}(\beta_0'x)>c}&=&\mathrm{O}_P(
\varepsilon_n),
\\
\sup_{x\in\mathcal{X}} \bigl|\nabla_{\beta}\hat{f} \bigl(\beta_0'x;
\beta_0 \bigr)-\nabla_{\beta}f \bigl(\beta_0'x;
\beta_0 \bigr)\bigr|\mathbf{1}_{f^{\tau
}_{\beta
_0}(\beta_0'x)>c} &=& \mathrm{O}_P
\bigl( \varepsilon'_n \bigr),
\end{eqnarray*}
%
where $\varepsilon_n$ and $\varepsilon_n'$ satisfy $\varepsilon_n\varepsilon'_n=\mathrm{o}(n^{-1/2})$,
$a_n^{-1/2}(\log n)^{1/2}\varepsilon_n \rightarrow0$ and
\mbox{$a_n^{-1/2}(\log n)^{1/2}\varepsilon'_n\rightarrow0$}.
\end{assum}



\section{Technical lemmas and proofs}\label{appB}

\renewcommand{\theequation}{A.\arabic{equation}}
\setcounter{equation}{5}

\renewcommand{\thelem}{A.\arabic{lem}}
\setcounter{lem}{0}

We start this Appendix with two technical lemmas, needed in the proofs
of the main results.
The first technical lemma gives a concentration inequality for the
convergence rate of semi-parametric estimators.

Let $b_n$ be a sequence of real numbers tending to zero, and let $\{
\zeta_{\alpha}\dvt \alpha\in\mathcal{A}\}$ be a family of uniformly
bounded functions, where $\mathcal{A}$ is a compact subset of $\mathbb
{R}^p$ (with $p \ge1$). Consider the class of functions
%
\begin{eqnarray}
\label{labelleclasse} && \mathcal{G} = \biggl\{(u,z,t,\delta) \mapsto
g_{\alpha
,x,v}(u,z,t,\delta)\nonumber
\\[-8pt]\\[-8pt]
 &&\hspace*{24pt}= K^0 \biggl(\frac{\psi(\alpha,u)-\psi(\alpha,x)}{b_n}
\biggr)\zeta_{\alpha
}(x,u,z,t,\delta)\xi(t)\mathbf{1}_{t\leq v} \dvt \alpha\in A, x \in\mathcal{X}, v \in\mathbb{R} \biggr\},\quad\nonumber
\end{eqnarray}
where $K^0$, $\psi$ and $\xi$ are fixed functions, $\mathcal
{X}\subset
\mathbb{R}^d$ is a compact set, and $t\in\mathbb{R},$ and consider
the process (in $\alpha$, $x$ and $v$)
\begin{eqnarray*}
\nu_n(g_{\alpha,x,v})=\sum_{i=1}^n
\bigl(g_{\alpha,x,v}(X_i,Z_i,T_i,
\delta_i)-E \bigl[g_{\alpha
,x,v}(X,Z,T,\delta ) \bigr] \bigr).
\end{eqnarray*}
Typically, $K^0$ denotes either a kernel or its derivative
of order 1 or 2. 

\begin{lem}
\label{talagrand} Assume that the class of functions
%
\begin{eqnarray}
\label{dernierehypothese} \biggl\{(u,z,t,\delta)\rightarrow K^0
\biggl(\frac{\psi(\alpha,u)-\psi(\alpha,x)}{b_n} \biggr)\zeta_{\alpha
}(x,u,z,t,\delta)\dvt \alpha
\in \mathcal{A}, x\in\mathcal{X} \biggr\}
\end{eqnarray}
is a VC-class of functions for a constant envelope, assume that
{ $E[|\xi(T)|^3]<\infty$}, and that {
$nb_n^3/(\log n)\rightarrow\infty$}. Then,
\begin{eqnarray*}
n^{-1/2}b_n^{-1/2}{ \bigl[\log(1/b_n)
\bigr]^{-1}} \|\nu_{n}\|_{\mathcal{G}}=
\mathrm{O}_P(1),
\end{eqnarray*}
where $\|\cdot\|_{\mathcal{G}}$ denotes the uniform norm over all maps
in $\mathcal{G}$.
\end{lem}

The proof of Lemma~\ref{talagrand} is a consequence of Proposition 1 in
Einmahl and Mason~\cite{r8} and Talagrand's inequality~\cite{r30}, and it is
available from the long version of this paper, see arXiv:\arxivurl{1111.6232}.

\begin{remark*} Note that if $K^0$ is of bounded variation with compact
support, and
if $\psi(\alpha,x)=\alpha'x,$ then (\ref{dernierehypothese}) holds, see
Nolan and Pollard~\cite{r24}.
\end{remark*}

The second technical lemma shows the consistency of the estimator $\hat
G_\theta(t|\lambda(\theta,x))$ and its vector of partial derivatives,
uniformly in $t,\theta$ and $x$, and it also establishes the rate of
convergence of the estimator $\hat{G}_{\hat{\theta}}(t|\hat{g}(x))$,
uniformly in $t$ and $x$.

\begin{lem}\label{uniftheta}
Under the assumptions of Theorem~\ref{curse3}, we have
%
\begin{eqnarray}
\label{eq1}\sup_{t\leq\tau, \theta\in\Theta, x\in
\mathcal{X}}\bigl|\hat{G}_{\theta} \bigl(t|\lambda(\theta,x)
\bigr)-G_{\theta
} \bigl(t|\lambda (\theta,x) \bigr)\bigr| &=&
\mathrm{o}_P(1),
\\
\label{eq2} \sup_{t\leq\tau, \theta\in\Theta, x\in
\mathcal{X}}\bigl|\nabla_{\theta}\hat{G}_{\theta}
\bigl(t| \lambda(\theta ,x) \bigr)-\nabla_{\theta}G_{\theta} \bigl(t|
\lambda( \theta,x) \bigr)\bigr| &=& \mathrm{o}_P(1),
\\
\label{eq3}\sup_{t\leq\tau}\sup_{x:g(x) \in
\mathcal{Z}_{\theta_0,\eta}}\bigl|\hat{G}_{\hat{\theta}}
\bigl(t|\hat {g}(x) \bigr)-G_{\theta_0} \bigl(t|g(x) \bigr)\bigr| &=&
\mathrm{O}_P \bigl(n^{-1/2}a_n^{-1/2}(
\log n)^{1/2} \bigr).\qquad
\end{eqnarray}
\end{lem}

\begin{pf}
For the first part, with probability tending to $1,$ for $t\leq
\tau,$ $1-\hat{G}(t|\lambda(\theta,x))>0.$ Taking the logarithm,
one obtains
\begin{eqnarray*}
\log \bigl(1-\hat{G} \bigl(t|\lambda(\theta,x) \bigr) \bigr) &=& \sum
_{i=1}^n (1-\delta_i)
\mathbf{1}_{T_i\leq t} \log \bigl(1-W_{n,i}(x,\theta) \bigr),
\end{eqnarray*}
where
\begin{eqnarray*}
W_{n,i}(x,\theta)&=&W_n(X_i,T_i;x,
\theta)\\
&=&{K \biggl(\frac{\lambda
(\theta,X_i)-\lambda(\theta,x)}{a_n} \biggr)}\bigg/\Biggl({\sum
_{j=1}^n \mathbf {1}_{T_j\geq T_i}K \biggl(
\frac{\lambda(\theta,X_j)-\lambda(\theta
,x)}{a_n} \biggr)}\Biggr).
\end{eqnarray*}
A Taylor expansion leads to
\begin{eqnarray*}
\log \bigl(1-\hat{G} \bigl(t|\lambda(\theta,x) \bigr) \bigr) &=& -\sum
_{i=1}^n (1-\delta_i)W_{n,i}(x,
\theta)\mathbf{1}_{T_i\leq
t}+\mathrm{O}_P \bigl(n^{-1}a_n^{-2}
\bigr),
\end{eqnarray*}
where the order of the remainder term is uniform in $t, \theta,
x,$ as
\[
\sup_{i:T_i\leq\tau}\sup_{x,\theta}\bigl|W_{n,i}(x,\theta )\bigr|=
\mathrm{O}_{P} \bigl(n^{-1}a_n^{-1}
\bigr).
\]
The remainder term is $\mathrm{o}_P(1)$
if $na_n^2\rightarrow\infty$. Rewrite
\begin{eqnarray*}
\sum_{i=1}^n (1-\delta_i)
\mathbf{1}_{T_i\leq t}W_{n,i}(x,\theta) &=&\frac{1}{na_n}\sum
_{i=1}^n(1-\delta_i)
\mathbf{1}_{T_i\leq t} K \biggl(\frac{\lambda(\theta,X_i)-\lambda(\theta,x)}{a_n} \biggr)S_{\theta
}
\bigl(\lambda(\theta,x),T_i \bigr)^{-1}
\\
&& {}+\frac{1}{na_n}\sum_{i=1}^n(1-
\delta_i)\mathbf{1}_{T_i\leq
t}K \biggl(\frac{\lambda(\theta,X_i)-\lambda(\theta,x)}{a_n} \biggr)
\\
&&\hspace*{45pt} {}\times\frac{\hat{S}_{\theta}(\lambda(\theta
,x),T_i)-S_{\theta}(\lambda(\theta,x),T_i)}{S_{\theta}(\lambda
(\theta
,x),T_i)\hat{S}_{\theta}(\lambda(\theta,x),T_i)},
\end{eqnarray*}
where
\begin{eqnarray*}
S_{\theta} \bigl(\lambda(\theta,x),y \bigr) &=& \bigl[1-H_{\theta}
\bigl(y|\lambda(\theta ,x) \bigr) \bigr]f_{Z_\theta} \bigl(\lambda(\theta,x)
\bigr),
\\
\hat{S}_{\theta} \bigl(\lambda(\theta,x),y \bigr) &=& \frac{1}{na_n}
\sum_{j=1}^n \mathbf{1}_{T_j\geq
y}K
\biggl(\frac{\lambda(\theta,X_j)-\lambda(\theta,x)}{a_n} \biggr).
\end{eqnarray*}
Apply Lemma~\ref{talagrand} to obtain the uniform
convergence of $\hat{S}_{\theta}$ towards $S_{\theta},$ and to show that
\begin{eqnarray*}
&&\sup_{x,\theta\in\Theta,t\leq
\tau} \Biggl\llvert \frac{1}{na_n}\sum
_{i=1}^n(1-\delta_i)\mathbf
{1}_{T_i\leq
t} K \biggl(\frac{\lambda(\theta,X_i)-\lambda(\theta,x)}{a_n} \biggr)S_{\theta
} \bigl(
\lambda(\theta,x),T_i \bigr)^{-1}
\\
&&\hspace*{41pt}-\int_{-\infty}^t \frac{\mathrm{d}H_{\theta,0}(s|\lambda(\theta,x))}{1-H_{\theta}(s-|\lambda
(\theta
,x))}
\Biggr\rrvert =\mathrm{o}_P(1).
\end{eqnarray*}
Since $S_{\theta}$ is uniformly bounded away from zero for $y\leq\tau,$
see Assumption~\ref{a1}, 
the result follows from
\[
\exp \biggl[-\int_{-\infty}^t\frac{\mathrm{d}H_{\theta,0}(s|\lambda(\theta
,x))}{1-H_{\theta}(s-|\lambda(\theta,x))}
\biggr] =1-G_{\theta} \bigl(t|\lambda(\theta,x) \bigr).
\]

For the gradient, we have
\begin{eqnarray*}
\nabla_{\theta}\hat{G}_{\theta} \bigl(t|\lambda(\theta,x) \bigr) &=&
\bigl(1-\hat {G}_{\theta} \bigl(t|\lambda(\theta,x) \bigr) \bigr)\sum
_{i=1}^n (1-\delta_i)
\mathbf{1}_{T_i\leq
t}\frac{\nabla_{\theta}W_{n,i}(x,\theta)}{1-W_{n,i}(x,\theta)}.
\end{eqnarray*}
From this, we deduce that the convergence of
$\nabla_{\theta}\hat{G}_{\theta}$ follows from the convergence of
$\hat
{G}_{\theta},$
of $\hat{S}_{\theta}$ and of
\begin{eqnarray*}
&& \frac{1}{na_n^2}\sum_{i=1}^n (1-
\delta_i)\mathbf{1}_{T_i\leq
t}\nabla_{\theta}\lambda(
\theta,x)K' \biggl(\frac{\lambda(\theta
,X_i)-\lambda(\theta,x)}{a_n} \biggr),
\end{eqnarray*}
and
\begin{eqnarray*}
&& \frac{1}{na_n^2}\sum_{i=1}^n
\mathbf{1}_{T_i\leq
t}\nabla_{\theta}\lambda(\theta,x)K'
\biggl(\frac{\lambda(\theta
,X_i)-\lambda(\theta,x)}{a_n} \biggr).
\end{eqnarray*}
These two quantities can be studied using Lemma~\ref{talagrand}, which
shows that their centered versions converge uniformly with rate
$(na_n^3)^{-1/2}\log n,$ while the bias term is of order $a_n^2.$

The third result can be deduced from a Taylor expansion,
Assumption~\ref{a13} and Proposition 4.3 in Van Keilegom and Akritas~\cite{r33}. Indeed, we can deduce that
\begin{eqnarray*}
&& \sup_{t\leq\tau}\sup_{x:g(x) \in
\mathcal{Z}_{\theta_0,\eta}}\bigl|\hat{G}_{\hat{\theta}} \bigl(t|\hat
{g}(x) \bigr)-G_{\theta_0} \bigl(t|g(x) \bigr)\bigr|
\\
&& \quad\leq \sup_{t\leq\tau}\sup_{x:g(x) \in
\mathcal{Z}_{\theta_0,\eta}}\bigl|\hat{G}_{\theta_0}
\bigl(t|g(x) \bigr)- G_{\theta_0} \bigl(t|g(x) \bigr)\bigr|+\mathrm{O}_P\bigl(
\|\hat{ \theta}-\theta_0\|\bigr).
\end{eqnarray*}
\upqed\end{pf}

We are now ready to give the proofs of the main results.

\begin{pf*}{Proof of Theorem~\ref{curse3}}
Part (i) of the theorem can be easily derived by replacing the
differentiability condition in Assumption~\ref{a13} by a uniform
continuity condition on $G_{\theta}$ with respect to $\theta,$ and
equation (\ref{eq1}) in Lemma~\ref{uniftheta}.

For part (ii), a Taylor expansion with
respect to $\theta$ leads to
\begin{eqnarray*}
\int \phi(x,y)\,\mathrm{d}[\hat{F}_{\hat{g}}-\hat{F}_g](x,y)=-
\frac{1}{n}\sum_{i=1}^n
\frac{\delta_i\phi(X_i,T_i)\nabla_{\theta}\hat{G}_{\theta
_n}(T_i-|\lambda(\theta_n,X_i))(\hat{\theta}-\theta_0)}{[1-\hat
{G}_{\theta_n}(T_i-|\lambda(\theta_n,X_i))]^2}
\end{eqnarray*}
for some $\theta_n$ between $\hat{\theta}$ and $\theta_0.$ From
the convergence of $\hat{\theta}$ towards $\theta_0,$ it follows that
$\theta_n$
tends to $\theta_0.$ Moreover, applying equation (\ref{eq1}) and
(\ref
{eq2}) in Lemma~\ref{uniftheta}, we obtain that
\begin{eqnarray*}
\int\phi(x,y)\,\mathrm{d}[\hat{F}_{\hat{g}}-\hat{F}_g](x,y)
&=& - \frac{1}{n}\sum_{i=1}^n
\frac{\delta_i\phi(X_i,T_i)\nabla_{\theta}G_{\theta
_0}(T_i-|\lambda
(\theta_0,X_i))(\hat{\theta}-\theta_0)}{[1-G_{\theta_0}(T_i-|g(X_i))]^2}
\\
&&{}+R_n(\phi),
\\
&=& U_n(\phi)+R_n(\phi),
\end{eqnarray*}
with $\sup_{\phi}|R_n(\phi)|\leq|R_n(\Phi)|=\mathrm{o}_P(n^{-1/2}),$ and
\begin{eqnarray*}
U_n(\phi) = \Biggl\{-\frac{1}{n}\sum
_{i=1}^n \frac{\delta_i\phi(X_i,T_i)\nabla_{\theta}G_{\theta
_0}(T_i-|\lambda
(\theta_0,X_i))}{[1-G_{\theta_0}(T_i-|g(X_i))]^2} \Biggr\} \Biggl\{
\frac
{1}{n}\sum_{j=1}^n
\mu(T_j,\delta_j,X_j) \Biggr
\}+R'_n(\phi),
\end{eqnarray*}
with $\sup_{\phi}|R'_n(\phi)|\leq|R'_n(\Phi)|=\mathrm{o}_P(n^{-1/2}).$
Centering the first sum in $U_n(\phi)$ and applying a uniform central
limit theorem (see, e.g., Van der Vaart and Wellner~\cite{r32}), we obtain
the stated representation.
\end{pf*}

\begin{pf*}{Proof of Theorem~\ref{th3}}
Consider the difference
\begin{eqnarray*}
&& \bigl|M_n(\beta,\hat{f},\tilde{J})-M_n(\beta,f,\tilde J)\bigr|
\\
&& \quad\leq2\int |y|\mathbf{1}_{y\leq\tau}\,\mathrm{d}\hat{F}_{\hat{g}}(x,y)
\sup_{x:\tilde
J(x)=1,\beta\in
\mathcal{B}}\bigl|\hat{f} \bigl(\beta'x;\beta \bigr)-f \bigl(
\beta'x;\beta \bigr)\bigr|
\\
&&\qquad{} +\int\mathbf{1}_{y\leq\tau}\bigl|\hat{f} \bigl(\beta'x;
\beta \bigr)+f \bigl(\beta'x;\beta \bigr)\bigr|\,\mathrm{d}
\hat{F}_{\hat{g}}(x,y) \sup_{x:\tilde
{J}(x)=1,\beta
\in
\mathcal{B}}\bigl|\hat{f} \bigl(
\beta'x;\beta \bigr)-f \bigl( \beta'x;\beta \bigr)\bigr|.
\end{eqnarray*}
%
%
The first term on the right-hand side converges uniformly to
zero by Assumption~\ref{premierevitesse} and the law of large numbers
for $\hat{F}_{\hat{g}}$ (see Theorem~\ref{t1} and Theorem~\ref{curse3}). The integral in the second term can be bounded by
\[
\bigl(1+\mathrm{o}_P(1) \bigr)\times\int2\Psi(x)\,\mathrm{d}
\hat{F}_{\hat{g}}(x,y),
\]
where $\mathrm{o}_P(1)$ is uniform in $\beta,$ by Assumption~\ref{etlemodele}
and~\ref{premierevitesse} -- (\ref{reg1}).
Now we have to show that $M_n(\beta,f,J^*)$ converges to
$M(\beta,f,J^*)$ uniformly in $\beta.$ For this, apply Theorem~\ref{t1}
and Theorem~\ref{curse3} using Assumption~\ref{etlemodele}. By usual arguments for
proving consistency (see, e.g., Van der Vaart~\cite{r31}, Theorem 5.7), the
consistency of $\beta_n$ follows.
\end{pf*}

\begin{pf*}{Proof of Lemma~\ref{esperanceestnulle}}
The proof is somewhat similar to the proof of Lemma 5A in Dominitz and Sherman~\cite{r6}.
First, observe that
\begin{eqnarray*}
f \bigl(\beta'X;\beta \bigr) = E \bigl[Y\mid
\beta'X,Y \leq \tau \bigr]
= E \bigl[f \bigl(\beta_0'X;\beta_0
\bigr)|\beta'X,Y\leq \tau \bigr]
=\frac{E [f(\beta_0'X;\beta_0)\mathbf{1}_{Y\leq\tau}\mid
\beta'X ]}{\mathbb{P}(Y\leq
\tau|\beta'X)}.
\end{eqnarray*}
Let $\alpha(X,\beta)=\beta_0'X-\beta'X$. Define
\[
\Gamma_X(\beta_1,\beta_2)=E \bigl[f \bigl(
\alpha(X,\beta_1)+\beta_2'X;
\beta_0 \bigr)\mathbf{1}_{Y\leq\tau}|\beta_2'X
\bigr],
\]
and note that $f(\beta'X;\beta) = \Gamma_X(\beta,\beta)/\mathbb{P}(Y
\le\tau|\beta'X)$.
Then,
\begin{eqnarray*}
\nabla_{\beta_1}\Gamma_X(\beta_0,
\beta_0) &= & -f' \bigl(\beta_0'X;
\beta_0 \bigr)E \bigl[X\mathbb{P} (Y\leq\tau\mid X )\mid
\beta_0'X \bigr],
\\
\nabla_{\beta_2}\Gamma_X(\beta_0,
\beta_0) &=& f' \bigl(\beta_0'X;
\beta_0 \bigr)X\mathbb{P} \bigl(Y\leq\tau\mid\beta_0'X
\bigr)+f \bigl(\beta_0'X;\beta_0 \bigr)
\nabla_{\beta} h(X,\beta_0),
\end{eqnarray*}
%
where $h(x,\beta)=\mathbb{P}(Y\leq\tau
|\beta'X=\beta'x).$ It follows that
%
\begin{eqnarray}
\label{tridule} &&\nabla_{\beta}f \bigl(\beta_0'x;
\beta_0 \bigr)
\nonumber
\\
&&\quad= \frac{f'(\beta_0'x;\beta_0) \{x\mathbb{P} (Y\leq\tau
\mid
\beta_0'X=\beta_0'x )-E [X\mathbb{P}(Y\leq\tau|X)\mid
\beta_0'X=\beta_0'x ] \}}{\mathbb{P}(Y\leq\tau
|\beta_0'X=\beta_0'x)}
\nonumber
\\[-8pt]
\\[-8pt]
&& \qquad {}+\frac{\nabla_{\beta}h(x,\beta_0)f(\beta_0'x;\beta_0)}{\mathbb
{P}(Y\leq
\tau|\beta_0'X=\beta_0'x)} -\frac{\nabla_{\beta}h(x,\beta_0)f(\beta_0'x;\beta_0)E[\mathbf{1}_{Y\leq\tau}|\beta_0'X=\beta_0'x]}{\mathbb
{P}(Y\leq
\tau|\beta_0'X=\beta_0'x)^2}
\nonumber
\quad
\\
&&\quad:= {xm_1 \bigl(\beta_0'x
\bigr)+m_2 \bigl(\beta_0'x \bigr).}
\nonumber
\end{eqnarray}
Therefore,
\begin{eqnarray*}
&& E \bigl[\nabla_{\beta}f \bigl(\beta_0'X;
\beta_0 \bigr)\mathbf{1}_{Y\leq
\tau
}|\beta_0'X
\bigr]
\\
&&\quad = \frac{E [f'(\beta_0'X;\beta_0) \{X\mathbb{P}
(Y\leq
\tau\mid
\beta_0'X )-E [X\mathbb{P}(Y\leq\tau|X)\mid
\beta_0'X ] \}\mathbf{1}_{Y\leq\tau} \mid\beta_0'X
]}{\mathbb{P}(Y\leq\tau|\beta_0'X)}
\\
&&\quad = 0.
\end{eqnarray*}
\upqed\end{pf*}

\begin{pf*}{Proof of Theorem~\ref{th4}} The proof consists of three steps:

%

\textit{Step 0: Replace $J$ by $J_0$.} { For
any $\mathcal{B}_n$ a sequence of shrinking neighborhoods of
$\beta_0$,
\[
\sup_{\beta\in\mathcal{B}_n} \bigl\llvert M_n(\beta,\hat f, J) -
M_n(\beta,\hat f, J_0) \bigr\rrvert \leq
\mathrm{o}_P \bigl(M_n(\beta,\hat f, J_0) +
n^{-1} \bigr).
\]
See Delecroix, Hristache and Patilea~\cite{r5}, page 738. Similar
arguments apply also when the trimming $J$ is defined with $\hat
f^\tau_{\beta_n} ( \beta'_n x)$ justifying the practical
implementation of the trimming function. }

\textit{Step 1: Bring the problem back to the parametric case.}

For notational simplicity, we work with $\nabla_\beta f$ instead of
$\nabla_{\tilde\beta} f$. Note that $\nabla_\beta f = (0,\nabla_{\tilde
\beta}'f)'$. We will show that, on $\mathcal{B}_n,$
\[
M_n(\beta,\hat{f},J_0)=M_n(
\beta,f,J_0)+\mathrm{o}_P \biggl(\frac{\|\beta
-\beta_0\|
}{\sqrt{n}}
\biggr)+\mathrm{o}_P \bigl(\|\beta-\beta_0
\|^2 \bigr)+C'_n,
\]
where $C_n'$ does not depend on $\beta.$ Decompose
\begin{eqnarray*}
M_{n} ( \beta,\hat{f},J_0 ) &=&M_{n} (
\beta,f,J_0 )
\\
&&{}-\frac{2}{n}\sum_{i=1}^{n}
\frac{\delta_{i}J_0(X_i) (
T_{i}-f ( \beta^{\prime}X_{i};\beta )  )\mathbf{1}_{T_i\leq\tau}
}{1-\hat
{G}_{\hat\theta} ( T_{i}-|\hat{g}(X_i) ) }%
\bigl[ \hat{f} \bigl( \beta^{\prime}X_{i};
\beta \bigr) -f \bigl( \beta^{\prime}X_{i};\beta \bigr) \bigr]
\\
&&{}+\frac{1}{n}\sum_{i=1}^{n}
\frac{\delta_{i}J_0(X_i)\mathbf{1}_{T_i\leq\tau}}{1-\hat{G}_{\hat\theta} (
T_{i}-|\hat{g}(X_i) ) } \bigl[ \hat{f} \bigl( \beta^{\prime
}X_{i};
\beta \bigr) -f \bigl( \beta^{\prime}X_{i};\beta \bigr)
\bigr]^{2}
\\
&=&M_{n} ( \beta,f,J_0 ) -2A_{1n}+B_{1n}.
\end{eqnarray*}

\textit{Step 1.1: Study of $A_{1n}.$}

$A_{1n}$ can be expressed as
\begin{eqnarray*}
A_{1n} &=& \frac{1}{n}\sum_{i=1}^{n}
\frac{\delta_{i}J_0(X_i) ( T_{i}-f ( \beta_{0}^{\prime
}X_{i};\beta_{0} )  ) \mathbf{1}_{T_i\leq
\tau}}{1-\hat{G}_{\hat\theta} ( T_{i}-|\hat{g}(X_i) )
} \bigl[ \hat{f} \bigl( \beta_{0}^{\prime}X_{i};
\beta_{0} \bigr) -f \bigl( \beta_{0}^{\prime}X_{i};
\beta_{0} \bigr) \bigr]
\\
&&{} +\frac{1}{n}\sum_{i=1}^{n}
\frac{\delta_{i}J_0(X_i)\mathbf
{1}_{T_i\leq\tau} (
f ( \beta_{0}^{\prime}X_{i};\beta_{0} ) -f ( \beta^{\prime}X_{i};\beta )  )
}{1-\hat{G
}_{\hat\theta} ( T_{i}-|\hat{g}(X_i) ) }
\\
&& \hspace*{34pt} {}\times \bigl[ \hat{f} \bigl( \beta^{\prime
}X_{i};
\beta \bigr) -f \bigl( \beta^{\prime}X_{i};\beta \bigr) -\hat{f}
\bigl( \beta_{0}^{\prime}X_{i};
\beta_{0} \bigr) +f \bigl( \beta_{0}^{\prime
}X_{i};
\beta_{0} \bigr) \bigr]
\\
&&{}+\frac{1}{n}\sum_{i=1}^{n}
\frac{\delta_{i}J_0(X_i)\mathbf
{1}_{T_i\leq\tau} (
f ( \beta_{0}^{\prime}X_{i};\beta_{0} ) -f ( \beta^{\prime}X_{i};\beta )  )
}{1-\hat{G%
}_{\hat\theta} ( T_{i}-|\hat{g}(X_i) ) } \bigl[ \hat {f} \bigl( \beta_{0}^{\prime}X_{i};
\beta_{0} \bigr) -f \bigl( \beta_{0}^{\prime}X_{i};
\beta_{0} \bigr) \bigr]
\\
&&{}+\frac{1}{n}\sum_{i=1}^{n}
\frac{\delta_{i}J_0(X_i) (
T_{i}-f
( \beta_{0}^{\prime}X_{i};\beta_{0} )  ) \mathbf{1}_{T_i\leq
\tau}}{1-\hat{G}_{\hat\theta} ( T_{i}-|\hat{g}(X_i) )
}
\\
&&\hspace*{34pt} {}\times \bigl[ \hat{f} \bigl( \beta^{\prime
}X_{i};
\beta \bigr) -f \bigl( \beta^{\prime}X_{i};\beta \bigr) -\hat{f}
\bigl( \beta_{0}^{\prime}X_{i};
\beta_{0} \bigr) +f \bigl( \beta_{0}^{\prime}X_{i};
\beta_{0} \bigr) \bigr]
\\
&=&A_{2n}+A_{3n}+A_{4n}+A_{5n}.
\end{eqnarray*}

$A_{2n}$ does not depend on $\beta.$ For $A_{3n},$ observe that,
for any $\beta\in\mathcal{B}_n,$ we can replace $J_0(X_i)$
by $\mathbf{1}_{f_{\beta}^{\tau}(\beta'X_i)>c/2}$ using
Assumption~\ref{trimberan}. As $\nabla_{\beta}f(\beta'x;\beta)$ is a
bounded function of $x$ and $\beta$ (Assumption~\ref{encoremodel},
since the class of functions has a bounded envelope), and using the
uniform convergence of $\nabla_{\beta}\hat{f}(\beta'x;\beta)$
(Assumption~\ref{premierevitesse}), we can obtain from a first order Taylor
expansion applied twice (for $f(\beta' x ;\beta)$ and for $\hat
f(\beta' x ;\beta) - f(\beta' x ;\beta)$ around $\beta_0$), that
$A_{3n}=\mathrm{o}_P(\|\beta-\beta_0\|^2).$

For $A_{4n},$ first replace $\hat{G}_{\hat\theta}$ with $G_{\theta_0}.$
For this, note that $[1-G_{\theta_0}(T_i-|g(X_i))]$ is bounded away
from zero with
probability tending to 1 for $T_i\leq\tau,$ and that
%
\begin{equation}
\label{dfr} \sup_{t\leq\tau,x:J_0(x)=1} \bigl\llvert \hat{G}_{\hat\theta} \bigl(t|
\hat {g}(x) \bigr)-G_{\theta_0} \bigl(t|g(x) \bigr) \bigr\rrvert \bigl\llvert
\hat{f} \bigl(\beta_0'x;\beta_0 \bigr)-f
\bigl( \beta_0'x;\beta_0 \bigr) \bigr\rrvert
=\mathrm{o}_P \bigl(n^{-1/2} \bigr)
\end{equation}
using part 2 of Assumption~\ref{unpeudevitesse?}, and Lemma~\ref{uniftheta}.
A first order Taylor expansion for $f(\beta'x;\beta)-f(\beta_0'x;\beta_0)$ and property (\ref{dfr}) lead to
\begin{eqnarray*}
A_{4n}&=& \frac{1}{n}\sum_{i=1}^n
\frac{\delta_{i}J_0(X_i)\mathbf{1}_{T_i\leq\tau} ( f ( \beta_{0}^{\prime}X_{i};\beta_{0} ) -f ( \beta^{\prime}X_{i};\beta )  )
}{1-G_{\theta_0}
( T_{i}-|g(X_i) ) } \bigl[ \hat{f} \bigl( \beta_{0}^{\prime}X_{i};
\beta_{0} \bigr) -f \bigl( \beta_{0}^{\prime
}X_{i};
\beta_{0} \bigr) \bigr]
\\
&& {}+\mathrm{o}_P \biggl(\frac{\|\beta-\beta_0\|}{\sqrt{n}} \biggr).
\end{eqnarray*}
Next, a second order Taylor development shows that
the first term above can be rewritten as
%
\begin{eqnarray}
\label{t} && \frac{(\beta-\beta_0)'}{n}\sum_{i=1}^n
\frac{\delta_{i}J_0(X_i)\mathbf{1}_{T_i\leq\tau}\nabla_{\beta}f(\beta_0'
X_i;\beta_0) }{1-G_{\theta_0}
( T_{i}-|g(X_i) ) } \bigl[ f \bigl( \beta_{0}^{\prime}X_{i};
\beta_{0} \bigr) -\hat{f} \bigl( \beta_{0}^{\prime
}X_{i};
\beta_{0} \bigr) \bigr]
\nonumber
\\[-8pt]
\\[-8pt]
&&\quad {}+\mathrm{o}_P \bigl(\|\beta-\beta_0
\|^2 \bigr).
\nonumber
\end{eqnarray}
To show that
this term is negligible, we will use empirical process theory.
We have that $f\in
\mathcal{H}_1,$ where $\mathcal{H}_1$ is the Donsker class
defined in Assumption~\ref{classededonsker}, and
$\hat{f}\in\mathcal{H}_1$ with probability tending to 1.
Consequently, the class of functions
\[
\mathcal{H}'_1= \biggl\{(y,c,x,t)\rightarrow
\frac{\mathbf{1}_{y\leq
c}\mathbf{1}_{y\leq\tau}\nabla_{\beta}f(\beta_0' x;
\beta_0)J_0(t)\phi(\beta_0't)}{1-G_{\theta_0} (y\wedge c-|g(x))} \dvt \phi\in\mathcal{H}_1 \biggr\}
\]
is a Donsker class, see Example
2.10.8 in Van der Vaart and Wellner~\cite{r32}. Furthermore, for all
$\phi\in\mathcal{H}_1,$
%
\begin{eqnarray}
\label{rouge}E \biggl[\frac{\delta
J_0(X)\nabla_{\beta}f(\beta_0'X;\beta_0)\phi(\beta_0'X)\mathbf
{1}_{T\leq
\tau}}{1-G_{\theta_0}
(T-|g(X))} \biggr]&=&E \bigl[\nabla_{\beta}f
\bigl(\beta_0'X;\beta_0 \bigr)\phi \bigl(
\beta_0'X \bigr)J_0(X)\mathbf{1}_{Y\leq
\tau}
\bigr]
\nonumber
\qquad
\\[-8pt]
\\[-8pt]
&=&0,
\nonumber
\end{eqnarray}
since $E[\nabla_{\beta}f(\beta_0'X,\beta_0)\mathbf{1}_{Y\leq
\tau}|\beta_0'X]=0$ (see Lemma~\ref{esperanceestnulle}), and since
$J_0(X)$ is a function of $\beta_0'X$ alone. Deduce that, since
$\mathcal{H}'_1$ is a Donsker class, and since $\hat{f}$ tends
uniformly to $f,$ that the first term in (\ref{t}) is of order
$\mathrm{o}_P(\|\beta-\beta_0\|n^{-1/2}).$ See the asymptotic equicontinuity
of Donsker classes, cf. Van der Vaart and Wellner~\cite{r32}, Section
2.1.2.

For $A_{5n},$ apply a second order Taylor expansion. Using that
$\nabla^2_{\beta,\beta}f$ is bounded, and that $\nabla^2_{\beta
,\beta
}\hat{f}$
converges uniformly to $\nabla^2_{\beta,\beta}f,$ we obtain
\begin{eqnarray*}
A_{5n} &=& \frac{(\beta-\beta_0)'}{n}\sum_{i=1}^n
\frac{\delta_iJ_0(X_i)\mathbf
{1}_{T_i\leq\tau}(T_i-f(\beta_0'X_i;\beta_0))[\nabla_{\beta
}f(\beta_0'X_i;\beta_0)-\nabla_{\beta}\hat{f}(\beta_0'X_i;\beta_0)]}{1-\hat
{G}_{\hat\theta}(T_i-|\hat{g}(X_i))}
\\
&& {}+\mathrm{o}_P \bigl(\|\beta-\beta_0
\|^2 \bigr).
\end{eqnarray*}
Proceed as for $A_{4n}$ to replace $\hat{G}$ and $\hat{g}$ by $G$
and $g,$ using part 3 of Assumption~\ref{unpeudevitesse?}. The
same arguments as for $A_{4n}$ can then be used, but considering
instead the
Donsker class
\[
\mathcal{H}'_2= \biggl\{(y,c,x)\rightarrow
\frac{\mathbf{1}_{y\leq
c}J_0(x)\mathbf{1}_{y\leq\tau}(y-f(\beta_0'x;\beta_0))\phi
(x)}{1-G_{\theta_0}(y-|g(x))}\dvt \phi\in\mathcal{H}_2 \biggr\},
\]
where $\mathcal{H}_2$ is defined in Assumption~\ref{classededonsker},
and observing that, for any function $\phi,$
\begin{eqnarray*}
&& E \biggl[\frac{\delta
J_0(X)\phi(X) (Y-f(\beta_0'X;\beta_0) )\mathbf{1}_{T\leq
\tau}}{1-G_{\theta_0}(T-|g(X))} \biggr]
\\
&&\quad =E \bigl[E \bigl[ \bigl(Y-f \bigl(\beta_0'X;
\beta_0 \bigr) \bigr)\mathbf{1}_{Y\leq\tau
}\mid X
\bigr]J_0(X)\phi(X) \bigr]=0,
\end{eqnarray*}
by the definition of our regression model.
Deduce that $A_{5n}=\mathrm{o}_P(\|\beta-\beta_0\|n^{-1/2}+\|\beta-\beta_0\|^2)$.

\textit{Step 1.2: Study of $B_{1n}.$}

Rewrite $B_{1n}$ as
\begin{eqnarray*}
B_{1n} &=&\frac{1}{n}\sum_{i=1}^{n}
\frac{\delta_{i}J_0(X_i)\mathbf{1}_{T_i\leq\tau}}{1-\hat{G}_{\hat\theta} (
T_{i}-|\hat{g}(X_i) ) }
\\
&&\hspace*{21pt} {}\times \bigl[ \hat{f} \bigl( \beta^{\prime
}X_{i};
\beta \bigr) -f \bigl( \beta^{\prime}X_{i};\beta \bigr) -\hat{f}
\bigl( \beta_{0}^{\prime}X_{i};
\beta_{0} \bigr) +f \bigl( \beta_{0}^{\prime
}X_{i};
\beta_{0} \bigr) \bigr]^{2}
\\
&&{}+\frac{1}{n}\sum_{i=1}^n
\frac{\delta_{i}J_0(X_i)\mathbf
{1}_{T_i\leq
\tau}}{1-\hat{G}_{\hat\theta} ( T_{i}-|\hat{g}(X_i) )
} \bigl[ \hat{f} \bigl( \beta_{0}^{\prime}X_{i};
\beta_{0} \bigr) -f \bigl( \beta_{0}^{\prime
}X_{i};
\beta_{0} \bigr) \bigr]^2
\\
&&{}+\frac{2}{n}\sum_{i=1}^n
\frac{\delta_{i}J_0(X_i)\mathbf
{1}_{T_i\leq
\tau}}{1-\hat{G}_{\hat\theta} ( T_{i}-|\hat{g}(X_i) )
} \bigl[ \hat{f} \bigl( \beta_{0}^{\prime}X_{i};
\beta_{0} \bigr) -f \bigl( \beta_{0}^{\prime}X_{i};
\beta_{0} \bigr) \bigr]
\\
&&\hspace*{33pt} {}\times \bigl[ \hat{f} \bigl( \beta^{\prime
}X_{i};
\beta \bigr) -f \bigl( \beta^{\prime}X_{i};\beta \bigr) -\hat{f}
\bigl( \beta_{0}^{\prime}X_{i};
\beta_{0} \bigr) +f \bigl( \beta_{0}^{\prime
}X_{i};
\beta_{0} \bigr) \bigr]
\\
&=& B_{2n}+B_{3n}+2B_{4n}.
\end{eqnarray*}

Observe that,
for any $\beta\in\mathcal{B}_n,$ we can replace $J_0(X_i)$
by $\mathbf{1}_{f_{\beta}^{\tau}(\beta'X_i)>c/2}$ using
Assumption~\ref{encoremodel}. Next, by a Taylor expansion and
the uniform convergence of $\nabla_{\beta}\hat{f}$, we have that
$B_{2n}=\mathrm{o}_P(\|\beta-\beta_0\|^2)$. The term
$B_{3n}$ does not depend on $\beta.$ For $B_{4n},$ a second order
Taylor expansion leads to
\begin{eqnarray*}
B_{4n}&=& \frac{(\beta-\beta_0)'}{n}\sum_{i=1}^n
\frac{\delta_{i}J_0(X_i) [ \hat{f} ( \beta_{0}^{\prime}X_{i};\beta_{0} ) -f ( \beta_{0}^{\prime
}X_{i};\beta_{0} )  ]}{1-\hat{G}_{\hat{\theta}} (
T_{i}-|\hat{g}(X_i) )}
\\
&& \hspace*{54pt} {}\times \bigl[\nabla_{\beta}\hat{f} \bigl(
\beta_0'X_i; \beta_0 \bigr)-
\nabla_{\beta
}f \bigl(\beta_0'X_i;
\beta_0 \bigr) \bigr]+\mathrm{o}_P \bigl(\|\beta-
\beta_0 \|^2 \bigr).
\end{eqnarray*}
Replace $\hat{G}$ by $G$ and use Assumption~\ref{unpeudevitesse?}, part
1, to conclude.

\textit{Step 2: Study of $M_n(\beta,f,J_0).$}

Observe that,
on $\mathrm{o}_P(1)$-neighborhoods of $\beta_0,$ from a Taylor expansion,
\begin{eqnarray*}
&& M_n(\tilde\beta,f,J_0)-M_n(\tilde
\beta_0,f,J_0)
\\
&& \quad=(\tilde\beta-\tilde\beta_0)'
\nabla_{\tilde\beta}M_n( \tilde \beta_0,f,J_0)+(
\tilde\beta-\tilde \beta_0)'\nabla^2_{\tilde\beta
,\tilde\beta
}M_n(
\tilde\beta_{0},f,J_0) (\tilde\beta-\tilde
\beta_0) +\mathrm{o}_P \bigl(\| \tilde \beta-\tilde
\beta_0\|^2 \bigr),
\end{eqnarray*}
and apply Theorem 1 and 2 of Sherman~\cite{r27} to conclude.
\end{pf*}

\begin{pf*}{Proof of Proposition~\ref{tt}}
The uniform convergence results in Assumptions~\ref{premierevitesse}
and~\ref{unpeudevitesse?} can be deduced from studying the uniform
convergence rate of the numerator and the denominator in~(\ref{f*})
(and their derivatives) separately.
This is a consequence of Lemma~\ref{talagrand}. Since the other terms
can be studied in a similar way, we only consider the case of the
denominator and its derivatives in (\ref{f*}). In each case, the bias
part can be dealt with uniformly with classical kernel arguments, and
is of order $h^2.$ For the centered version of~$f^*,$ the result can be
deduced from the study of the uniform convergence rate of empirical
processes indexed by some class of functions as the one defined in
(\ref{labelleclasse}), with
\[
\zeta_{\beta}(x,X,Z,T,\delta)=\frac{\delta(x-X)^{j}}{1-G_{\theta
_0}(T-|Z)},
\]
where $j=0$ (resp., 1, 2) for $f^*$ (resp., $\nabla_{\beta}f^*$, $\nabla^2_{\beta,\beta}f^*$), and $\xi(T)=T.$ The kernel $K^0$ in (\ref
{labelleclasse}) is either $\tilde K$ or $\tilde K'$ or $\tilde K'',$
and $\psi(\beta,x)=\beta'x.$ It follows from the conditions on
$\tilde K$
and from Nolan and Pollard~\cite{r24} that the class of functions
\[
\biggl\{x\rightarrow K^0 \biggl(\frac{\beta'x-\beta'u}{{h}} \biggr)\dvt u\in
\mathcal{X}, {h>0,\beta\in\mathcal{B}} \biggr\}
\]
is a VC-class of functions. Moreover, $u\rightarrow(x-u)^j$
($j=0,1,2$) is also a VC-class of bounded functions
{using permanence properties of VC-classes, see Lemma
2.6.18 in Van der Vaart and Wellner~\cite{r32}.} Finally, since
$1-G_{\theta_0}(T-|Z)$ is bounded away from zero,
(\ref{dernierehypothese}) holds. Now applying Lemma~\ref{talagrand},
we get
\begin{eqnarray*}
\sup_{\beta,x}\bigl|f^* \bigl(\beta'x;\beta \bigr)-f \bigl(
\beta'x;\beta \bigr)\bigr| &=& \mathrm{O}_{P} \bigl((\log
n)^{1/2}n^{-1/2}h^{-1/2}+h^2 \bigr),
\\[-2pt]
\sup_{\beta,x}\bigl|\nabla_{\beta}f^* \bigl(\beta'x; \beta
\bigr)-\nabla_{\beta
}f \bigl(\beta'x;\beta \bigr)\bigr| &=&
\mathrm{O}_{P} \bigl((\log n)^{1/2}n^{-1/2}h^{-3/2}+h^2
\bigr),
\\[-2pt]
\sup_{\beta,x}\bigl|\nabla_{\beta,\beta}^2f^* \bigl(
\beta'x;\beta \bigr)-\nabla_{\beta
,\beta}^2f \bigl(
\beta'x;\beta \bigr)\bigr| &=& \mathrm{O}_{P} \bigl((\log
n)^{1/2}n^{-1/2}h^{-5/2}+h^2 \bigr),
\end{eqnarray*}
where $h^2$ comes from the bias term. Hence, Assumption
\ref{premierevitesse} holds if $h\rightarrow0$ and
$nh^5(\log n)^{-1/2}\rightarrow\infty.$ Assumption
\ref{unpeudevitesse?} holds if $(\log
n)^{-1}n^{1/2}a_n^{1/2}h\rightarrow\infty,$ and $nh^8\rightarrow
0$.

The first part of Assumption~\ref{classededonsker} follows directly
from the uniform convergence of $f^*.$ Elementary algebra shows that
the gradient of $f^*$ can be written as
%
\begin{equation}
\label{decompof*} \nabla_{\beta}f^* \bigl(\beta_0'x;
\beta_0 \bigr)=xm^*_1 \bigl(\beta_0'x
\bigr)+m^*_2 \bigl(\beta_0'x \bigr).
\end{equation}
Using the same arguments as above, these two functions converge
uniformly to $m_1(\beta_0'x)$ and $m_2(\beta_0'x)$, respectively,
where
$\nabla_{\beta}f(\beta_0'x;\beta_0)=xm_1(\beta_0'x)+m_2(\beta_0'x),$
{see equation (\ref{tridule})}, and Assumption
\ref{classededonsker} follows.
\end{pf*}

The proof of Proposition~\ref{diff} is a direct consequence of Lemma
\ref{uniftheta}, equation (\ref{eq3}), and the fact that
\[
\sup_{\beta} \Biggl[(nh)^{-1}\sum_{i=1}^n
\bigl|\tilde K^{(j)}\bigr| \biggl(\frac
{\beta'x-\beta'X_i}{h} \biggr)|T_i|^k
\Biggr] = \mathrm{O}_P(1),
\]
and hence will be omitted.
\end{appendix}
%

\section*{Acknowledgements}

The third author acknowledges financial support from IAP research
network P6/03 of the Belgian Government (Belgian Science Policy), and
from the European Research Council under the European Community's
Seventh Framework Programme (FP7/2007-2013)/ERC Grant agreement No.
203650.



\printhistory

\end{document}